\documentclass[11pt,reqno,twosdie]{amsart}

\usepackage{amsmath,amssymb,mathrsfs}
\usepackage{pifont}
\usepackage{floatrow}
\newfloatcommand{capbtabbox}{table}[][\FBwidth]

\usepackage[asymmetric,top=3.2cm,bottom=3.9cm,left=2cm,right=2cm]{geometry}
\usepackage{amsmath,amsfonts,amsthm,mathrsfs,amssymb,cite}
\usepackage[usenames]{color}
\usepackage{mathtools} %for \mathclap, \mathrlap, \smashoperator etc
\usepackage{graphics}
\usepackage{framed}
\usepackage{enumerate}
\usepackage{subfigure}
\usepackage{footmisc}
\usepackage{hyperref}
\usepackage{xcolor}
\usepackage{caption}
\hypersetup{
  colorlinks,
  linkcolor={blue!80!black},
  urlcolor={blue!80!black},
  citecolor={blue!80!black}
}

\usepackage{graphicx}
\usepackage{epstopdf}
\usepackage{latexsym}
\usepackage{enumerate}
\usepackage{pgfplots, tikz}
\pgfplotsset{width = 10 cm}
\theoremstyle{plain}
\newtheorem{theorem}{Theorem}[section]
\newtheorem{lemma}[theorem]{Lemma}

\newtheorem{remark}[theorem]{Remark}

\newtheorem{algo}[theorem]{Algorithm}
\renewcommand{\eqref}[1]{\textnormal{(\ref{#1})}}

\numberwithin{equation}{section}

\usepackage{color}
\usepackage{booktabs}

\renewcommand{\Im}{{\ensuremath{\mathrm{Im\,}}}}

\renewcommand{\Re}{{\ensuremath{\mathrm{Re\,}}}} %Realteil nicht als fraktur

\title[\resizebox{4.6in}{!}{On the detection of medium inhomogeneity by contrast agent}]{On the detection of medium inhomogeneity by contrast agent: wave scattering models and numerical implementations}

\author[Zhe Wang,  Ghandriche and Jijun Liu]
{Zhe Wang$^{\dag,*}$ Ahcene Ghandriche$^{*}$ and Jijun Liu$^{\dag,*,\ddag}$}

\thanks{$ ^{\dag}$School of Mathematics, Southeast University, Nanjing, 210096, P.R.China.}

\thanks{\hskip 0.35cm$ ^{*}$Nanjing Center for Applied Mathematics, Nanjing, 211135, P.R.China.}

\thanks{\hskip 0.35cm$ ^{\ddag}$Corresponding author. Prof. Dr. Jijun Liu}

\thanks{\hskip 0.35cm This work was supported by NSFC (No.12241102) and also the Jiangsu Provincial Scientific Research Center of Applied Mathematics under Grant No.BK20233002.}

\thanks{\hskip 0.35cm Email address: zhewang23@seu.edu.cn (Zhe Wang) gh.hsen.math@gmail.com (Ahcene Ghandriche)  jjliu@seu.edu.cn (Jijun Liu) }

\usepackage{amsmath,amssymb,mathrsfs}

\usepackage[asymmetric,top=3.2cm,bottom=3.9cm,left=2cm,right=2cm]{geometry}
\usepackage{amsmath,amsfonts,amsthm,mathrsfs,amssymb,cite}
\usepackage[usenames]{color}
\usepackage{mathtools} %for \mathclap, \mathrlap, \smashoperator etc
\usepackage{graphics}
\usepackage{framed}

\usepackage{enumerate}

\usepackage{subfigure}

\usepackage{hyperref}
\usepackage{xcolor}
\hypersetup{
  colorlinks,
  linkcolor={blue!80!black},
  urlcolor={blue!80!black},
  citecolor={blue!80!black}
}
\usepackage{graphicx}
\usepackage{latexsym}
\usepackage{enumerate}

\theoremstyle{plain}
\renewcommand{\eqref}[1]{\textnormal{(\ref{#1})}}

\numberwithin{equation}{section}

\usepackage{color}
\usepackage{booktabs}
\renewcommand{\Im}{{\ensuremath{\mathrm{Im\,}}}}

\renewcommand{\Re}{{\ensuremath{\mathrm{Re\,}}}} %Realteil nicht als fraktur

\date{\today}

\begin{document}

\begin{abstract}
We consider the wave scattering and inverse scattering in an inhomogeneous medium embedded a homogeneous droplet with a small size, which is modeled by a constant mass density and a small bulk modulus. Based on the Lippmann-Schwinger integral equation for scattering wave in inhomogeneous medium, we firstly develop an efficient approximate scheme for computing the scattered wave as well as its far-field pattern for any droplet located in the inhomogeneous background medium.
By establishing the approximate relation between the far-field patterns of the scattered wave before and after the injection of a droplet, the scattered wave of the inhomogeneous medium after injecting the droplet is represented by a measurable far-field patterns, and consequently the inhomogeneity of the medium can be reconstructed from the Helmholtz equation. Finally, the reconstruction process in terms of the dual reciprocity method is proposed to realize the numerical algorithm for recovering the bulk modulus function inside a bounded domain in three dimensional space, by moving the droplet inside the bounded domain. Numerical implementations are given using the simulation data of the far-field pattern to show the validity of the reconstruction scheme, based on the mollification scheme for dealing with the ill-posedness of this inverse problem.

\bigskip



\end{abstract}

\maketitle

\section{Introduction}
The wave propagation in inhomogeneous medium corresponding to given incident wave is an important physical phenomenon.  For different kinds of waves, the governed equation may be in different forms, say, the Helmholtz equation for acoustic wave, the Maxwell equations for electromagnetic wave, and the elastic wave equations for elastic waves. In general, the inhomogeneity of the medium is represented by some physical parameters such as medium density and wave velocity arising in the governed equations. The reconstruction of the medium parameters using some information about the scattered wave such as its far-field pattern, which are called the inverse medium problems, have been thoroughly studied, see \cite{colton2019inverse} and the references therein.  It has been proven in  \cite[Theorem 10.5]{colton2019inverse} that the far-field pattern of the scattered wave using all observation points excited by incident plane waves of the form
\begin{equation*}
    u^{inc}(x,\theta)  :=  e^{i  k  x \cdot \theta}  =  e^{i  k \left\vert x \right\vert  \hat{x} \cdot \theta}, \quad x \in \mathbb{R}^{3}, \;\; \theta, \hat{x} \in \mathbb{S}^2:=\{x\in \mathbb{R}^3: |x|=1\},
\end{equation*}
where $k$ is the fixed wave number, $\theta$ is the incident direction, and $\hat{x}$ is the unit observation direction, from all incident directions $\theta$ can uniquely determine the refraction index with compact support in $\mathbb{R}^3$. Of course, the requirement for the unique determination of the medium inhomogeneity that both $\hat x$ and $\theta$ should be full of the unit sphere $\mathbb{S}^2$ is not practical from the engineering point of view.

Except for the wave propagation in general inhomogeneous media, the other important physical situations are the wave propagation in inhomogeneous media embedded with some bubbles, or droplets. Many interesting physical works have been devoted to the acoustic bubble problems, see
\cite{DGS, HA, HABF}. In the case that the bubbles are of small size, they will yield extra scattering behaviors like the point-source.
Then by moving bubbles in the inhomogeneous medium, the corresponding scattering waves depending on the bubbles will involve information about the background media, which provide a new way to the inhomogeneity detection for the background media. That is, by applying the far-field patterns of the scattered waves in inhomogeneous background medium embedded with various bubbles, it is possible to recover the refraction index of the medium using finite number of incident directions and observation points in $\mathbb{S}^2$, instead of all incident directions and all observation directions. Such a configuration with finite number of incident directions and observations for the inverse medium problems based on the wave scattering is of crucial importance.

For $\rho_{0}(\cdot), k_{0}(\cdot)\in C(\mathbb{R}^3)$ representing the parameters of the inhomogeneity of the background medium, as the usual situations, we restrict the inhomogeneity  in a bounded domain $\Omega\subset\subset \mathbb{R}^3$. More precisely, we assume that
$\rho_0(x)=1$ and $k_0(x)=1$ for $x\in \mathbb{R}^3 \setminus  \Omega$,
where $\rho_{0}(\cdot)$ (respectively, $k_{0}(\cdot)$) is the density function (respectively, bulk modulus function) of the background medium.

Here we assume that $k_0(x)=1$, in $\mathbb{R}^3 \setminus  \Omega$, is just for the simplicity of the statement. In principle, such an assumption can be relaxed to the case that the inhomogeneity of $k_0(\cdot)$ is in the full space with general distribution. We contend that the analysis developed herein can be extended straightforwardly, leading to analogous results up to a multiplicative constant. In the practical situation that
$k_{0}(x)  =  1  +  \mathcal{O}\left( \left\vert x \right\vert^{-\alpha}\right)$ with some $\alpha > 0$,
which means that the medium inhomogeneity is restricted almost in a bounded domain, we can begin by embedding the imaging domain $\Omega$ in a large ball $B({\bf 0},R)$, and then use the same algorithm to recover $k_{0}(\cdot)$, in $B({\bf 0},R)$. This introduces an additional error term  given by  $\mathcal{O}(R^{-\alpha})$. However, this process is another important research topic.

The effective equations for wave propagation in bubbly liquids have been derived in the quasi-static regime, see \cite{CMPT, SGK} for low frequency regime. In this paper, however, we are more concerned with wave propagation in the resonant regime, i.e., moderate frequency regime, from a numerical perspective.
We consider the scattering of acoustic waves in  inhomogeneous media in $\mathbb{R}^{3}$ embedded a homogeneous droplet $D_z$ given by the form
\begin{equation}\label{DzepsB}
 D_z  :=  z  +  \varepsilon  B({\bf{0}},1),
\end{equation}
where $B({\bf{0}},1)$ is the unit ball with center ${\bf 0}$ in $\mathbb{R}^3$, with its boundary $\partial B({\bf{0}},1)\equiv \mathbb{S}^2$, and $z$ specifies the location of the droplet with radius $0 < \varepsilon \ll 1$, i.e., $D_z=B(z, \varepsilon)$. The small parameter $0  < \varepsilon \ll 1$ means that the used droplet is of small size, and consequently acts like a particle, or, a point source. Then the mass density and bulk modulus $(\rho(\cdot),k(\cdot))$ in $\mathbb{R}^3$ are of the form
\begin{equation}\label{defrhodefk}
       \rho(\cdot)  :=  \begin{cases}
			\rho_{0}(\cdot) & \text{in} \; \mathbb{R}^3 \setminus D_z\\
            \rho_{1} & \text{in} \; {\overline D_z}
		 \end{cases}
   \quad\text{and} \quad
   k(\cdot)  :=  \begin{cases}
			k_{0}(\cdot) & \text{in} \; \mathbb{R}^3 \setminus D_z\\
            k_{1} & \text{in} \; {\overline D_z}
		 \end{cases},
\end{equation}
where $\rho_{1}$ and $k_{1}$ are known positive constants  indicating the constant density and the bulk modulus of the homogeneous droplet $D_z$.

For a given incident plane wave
$u^{inc}(x, \theta, \omega)  := e^{i  \omega  x \cdot \theta}$, where $\omega$ is the frequency and $x \in \mathbb{R}^3$, solving
\begin{equation}\label{PDE-u-inc}
\Delta u^{inc}(x, \theta, \omega) +  \omega^{2}  u^{inc}(x, \theta, \omega)  =  0, \quad x \in \mathbb{R}^{3},
\end{equation}
the total wave $u(\cdot, \theta, \omega)  :=  u^{s}(\cdot, \theta, \omega)  + u^{inc}(\cdot, \theta, \omega)$, propagating in $\mathbb{R}^3$ embedded a droplet scatter $D_z$ with media parameters given by (\ref{defrhodefk}),
obeys the governed system
\begin{eqnarray}\label{Eq0758}
\begin{cases}
    \nabla \cdot \left[ \dfrac{1}{\rho_{0}} \nabla u \right]  + \dfrac{\omega^{2}}{k_{0}} u  =  0 & \text{in} \;  \mathbb{R}^3\setminus \overline{D}_z, \\
    \nabla \cdot \left[ \dfrac{1}{\rho_{1}} \nabla u \right]  +  \dfrac{\omega^{2}}{k_{1}} u  =  0 & \text{in} \;  D_z, \\
    u|_{-}  -  u|_{+}  =  0 & \text{on} \; \partial D_z,\\
    \dfrac{1}{\rho_{1}}  \dfrac{\partial u}{\partial \nu}|_{-}  -  \dfrac{1}{\rho_{0}}  \dfrac{\partial u}{\partial \nu}|_{+}  =  0 & \text{on} \;\; \partial D_z,
\end{cases}
\end{eqnarray}
where $\nu$ denotes the outward unit normal to $\partial D_z$,
and $u^{s}(\cdot,\cdot,\cdot)$ is the corresponding scattered wave satisfying the Sommerfield radiation condition at infinity, i.e.,
\begin{equation}\label{liu02}
    \frac{\partial u^{s}}{\partial \left\vert x \right\vert}(x, \theta, \omega)  -  i  \omega  \sqrt{\frac{ \rho_{0}(x)}{k_{0}(x)}}  u^{s}(x, \theta, \omega)  =  \mathcal{O} \left( \frac{1}{\left\vert x \right\vert} \right), \quad \left\vert x \right\vert \rightarrow + \infty.
\end{equation}
The problem \eqref{Eq0758}-\eqref{liu02} is well-posed, see \cite{ACP, ACP-2}. The scattered wave exhibits the asymptotic behavior
\begin{equation*}
    u^{s}(x, \theta, \omega)  =  \frac{e^{i  \omega  \sqrt{\frac{\rho_{0}(x)}{k_{0}(x)}}  \left\vert x \right\vert}}{\left\vert x \right\vert}  u^{\infty}(\hat{x}, \theta, \omega)  +  \mathcal{O}\left( \frac{1}{\left\vert x \right\vert^{2}} \right), \quad \left\vert x \right\vert \rightarrow + \infty,
\end{equation*}
where $\hat{x} := \frac{x}{\left\vert x \right\vert} \in \mathbb{S}^{2}$ and $u^{\infty}(\cdot, \theta, \omega)$ is the far-field pattern corresponding to $u^{s}(\cdot, \theta, \omega)$. Of course, both $u^{s}(\cdot,\theta, \omega)$ and $u^{\infty}(\cdot,\theta, \omega)$ depend on $D_z$, which will be denoted by $u^{s}_z(\cdot,\theta, \omega), u^{\infty}_z(\cdot,\theta, \omega)$ in the sequel. Moreover,
when we move $z\in\Omega$ and keep $D_z\subset\subset \Omega$, the scattered wave $\{u^s_z(x,\theta, \omega):\;  x\in \mathbb{R}^{3}, \theta \in \mathbb{S}^{2}\}$  involves the information about $k_0(x)$ in $\{x\in\Omega:\; \hbox{dist}(x,\partial\Omega)\ge \varepsilon>0\}$.

The main theme of this paper is the formalization of this extraction process as a numerical scheme. An important feature of this reconstruction process is that we can only recover $k_0(x)$ in the interior part of $\Omega$, since our numerical reconstruction scheme can be realized efficiently only in the domain $\{x\in\Omega:\; \hbox{dist}(x,\partial\Omega)\ge \varepsilon>0\}$. The reason is that our key formulas are valid only when the frequency is sufficiently close to the eigenvalues of the Newtonian potential operator $N_{D}(\cdot)$, see $(\ref{Def-N_B-Equa0708})$ for its definition, coming from the fact that  the eigenvalues of $N_{D}(\cdot)$ become unbounded when the domain $D$ approaches the outer boundary $\partial \Omega$  in spectral theory. Consequently, as we cannot approach unbounded eigenvalues using the same ideas, our proposed numerical algorithm fails in this regime. To circumvent this issue, we have assumed that $D$ is situated away from the boundary $\partial \Omega$.  The other important observation is that the scattered wave near $z$ for the perturbed medium with droplet $D_z$ is sensitive to $k_{0}(\cdot)$. An intuitive explanation on this phenomenon is that, since the droplet $D_z$ with small size $\varepsilon$ behaviors like a point source, the corresponding scattered wave in $z$ is theoretically infinity, and consequently very hard to be approximated.

So an important inverse problem is to recover $k_{0}(\cdot)$ inside $\Omega$ using the information about $u^s_z(\cdot,\theta, \omega)$ for given fixed incident direction $ \theta \in \mathbb{S}^{2}$ and all possible $D_z  \subset \subset \Omega$. By this configuration, for a known domain $\Omega$ where the background medium is inhomogeneous, we can detect its inhomogeneity.

In the absence of a droplet, \eqref{Eq0758} in the case $\rho_0(x)\equiv 1$ in $\mathbb{R}^3$ is reduced to the  Helmholtz equation
\begin{equation}\label{Heq}
    \Delta v(x,\theta, \omega)   +  \frac{\omega^{2}}{k_{0}(x)}  v(x,\theta, \omega)   =  0, \quad x \in \mathbb{R}^{3},
\end{equation}
where $v(\cdot,\theta, \omega)  =  v^s(\cdot,\theta, \omega)  +  u^{inc}(\cdot,\theta, \omega)$ is the total field before injecting the droplet $D_z$. The scattered wave $v^s(\cdot,\theta, \omega)$ is also of the asymptotic behavior
\begin{equation*}
    v^{s}(x, \theta, \omega)  =  \frac{e^{i  \omega  \sqrt{\frac{1}{k_{0}(x)}}  \left\vert x \right\vert}}{\left\vert x \right\vert}  v^{\infty}(\hat{x}, \theta, \omega)  +  \mathcal{O}\left( \frac{1}{\left\vert x \right\vert^{2}} \right), \quad \left\vert x \right\vert \rightarrow +\infty,
\end{equation*}
where $v^{\infty}(\cdot, \theta, \omega)$ is the far-field pattern corresponding to the scattered field $v^{s}(\cdot, \theta, \omega)$.

The inverse problem we are facing is the recovery of $k_{0}(\cdot)$ in $\Omega$ from the data $u_z^{\infty}(\cdot,\cdot,\cdot)$, corresponding to all $D_z\subset\subset\Omega$, and the data $v^\infty(\cdot,\cdot,\cdot)$. To handle this, thanks to \cite[section 2]{DGS}, the collection of data in a single backward scattering direction near a selected resonance was shown to be enough to recover the total field $v(\cdot,\cdot,\cdot)$, and therefore the bulk modulus function $k_{0}(\cdot)$. Our motivation to reduce the complexity of our problem from the computational and numerical perspective, see \textbf{Remark \ref{COEF}}, motivates us to follow the same idea suggested in \cite[section 2]{DGS} by applying one incident direction $\theta$ and one observation direction along $-\theta$, i.e., the back-scattered direction.

This paper is arranged as follows. In section \ref{SectionI-II}, we set up a reduced acoustic model corresponding to a constant mass density and variable bulk modulus. In section \ref{SectionII}, we will write done the computations related to the far-field patterns $v^{\infty}(\cdot,\cdot,\cdot)$ and $u_z^{\infty}(\cdot,\cdot,\cdot)$. In subsection \ref{SubSectionUnperturbed}, we explain how we can estimate the total field $v(\cdot,\cdot,\cdot)$ and the far-field $v^{\infty}(\cdot,\cdot,\cdot)$ for the unperturbed medium cases. In subsection \ref{SubSecIII-III}, we explain how we can estimate the total field $u_z(\cdot,\cdot,\cdot)$ and the far-field $u_z^{\infty}(\cdot,\cdot,\cdot)$ for the perturbed medium cases. The objective of section \ref{Sec-Num-Sim} is to numerically support the derived results in section \ref{SectionII}. Finally, section \ref{SectionV} is devoted to the numerical reconstruction of the bulk modulus function $k_{0}(\cdot)$ in interior domain of $\Omega$, where the mollification scheme is applied to deal with the ill-posedness of this inverse problem, due to the numerical differentiations in the reconstruction algorithm.

\section{Reduced Acoustic Model}\label{SectionI-II}
To recover $k_0(\cdot)$ in $\Omega$, we introduce a droplet scatter $D_z\subset\subset \Omega$ for $z\in\Omega$ and consider the wave scattering restricted in $\Omega$.  Under the assumption that the density function $\rho(\cdot)$ given by \eqref{defrhodefk} is constant, which means the inhomogeneity of the scattering media and the droplet is independent of the media density, the corresponding mathematical model \eqref{Eq0758} becomes
\begin{eqnarray}\label{Eq1}
	\begin{cases}
	\Delta  u_{z}  +  \frac{\omega^2}{k_0}  u_{z}  =  0 &\mbox{in}\ \Omega\setminus \overline{D}_z,\\
 	\Delta u_{z}  +  \frac{\omega^2}{k_1}  u_{z}  =  0 &\mbox{in} \;  D_z,\\
		u_{z}|_{-}  -  u_{z}|_{+}  =  0 &\mbox{on} \; \partial D_z, \\
      \frac{\partial u_{z}}{\partial_{\nu}}|_{-}
 		  -  \frac{\partial u_{z}}{\partial_{\nu}}|_{+} =0 \quad &\mbox{on} \; \partial D_z,
	\end{cases}
\end{eqnarray}
 where we have assumed
$\rho(\cdot) \equiv 1$ in $\mathbb{R}^{3}$.
The notation $u_{z}(\cdot,\cdot,\cdot)$ denotes the total field  after injecting a droplet $D_z$. To establish our reconstruction scheme,  we assume that the droplet is of a small bulk modulus of order $\varepsilon^{2}$, that is,
\begin{equation}\label{Ctek1}
k_{1} =  \overline{k}_{1}  \varepsilon^{2},
\end{equation}
where $\overline{k}_{1}$ is some constant independent of $\varepsilon$, and $\varepsilon$ is the radius of $D_z$.
The following results about the far field patterns of $u^{s}_{z}(\cdot, \cdot, \cdot)$ and $v^{s}(\cdot,\cdot,\cdot)$ have been established in \cite[Theorem 1.2]{DGS}.

\begin{lemma}
Denote by $\{ (\lambda_{n}, e_{n}),\;n \in \mathbb{N}\}$  the eigensystem related to the Newtonian operator $N_{B}$, defined from $\mathbb{L}^{2}(B)$ to $\mathbb{L}^{2}(B)$ by
\begin{equation}\label{Def-N_B-Equa0708}
    N_{B}(f)(x)  := \int_{B} \frac{1}{4  \pi \left\vert x - y \right\vert}  f(y)  dy, \quad x \in B:=B({\bf 0},1),
\end{equation}
 i.e., $N_{B}(e_{n})  =  \lambda_{n} e_{n}$ in $B$. Then, for $z\in\Omega$ such that $D_z\subset\subset\Omega$, the far-filed patterns admit
\begin{equation}\label{DiffUinfvinf}
u^{\infty}_z(\hat{x},\theta,\omega)  =  v^{\infty}(\hat{x},\theta,\omega)  -  \frac{1}{4 \pi  \overline{k}_{1}}  \frac{\omega^{2}  \omega_{n_{0}}^{2}}{\left( \omega^{2}  -  \omega_{n_{0}}^{2} \right)}  \left( \int_B e_{n_{0}}(x)  dx \right)^{2} \varepsilon  v(z,-\hat{x},\omega)  v(z,\theta, \omega)  +  \mathcal{O}\left( \varepsilon^{\min(1,2-2h)} \right)
\end{equation}
uniformly for all $\theta,\hat{x}\in \mathbb{S}^{2}$,  under the condition $\omega^{2}  -  \omega^{2}_{n_{0}}  =  \mathcal{O}\left( \varepsilon^{h} \right)$ for $h\in (0,1)$, where $\omega^{2}_{n_{0}}$ is the eigen-frequency given by
\begin{equation}\label{Def-Eigen-Freq}
    \omega^{2}_{n_{0}}  =  \overline{k}_{1}  \lambda_{n_0}^{-1}.
\end{equation}
\end{lemma}

Then, by taking the observation direction $\hat{x}$ as the back scattered direction, the equation \eqref{DiffUinfvinf} yields
\begin{equation}\label{BSD}
u_z^{\infty}(-  \theta, \theta, \omega)  =  v^{\infty}(-  \theta, \theta, \omega)  - \frac{1}{4  \pi \overline{k}_{1}}  \frac{\omega^{2}  \omega_{n_{0}}^{2}}{\left( \omega^{2} -  \omega_{n_{0}}^{2} \right)}  \left( \int_{B} e_{n_{0}}(x)  dx \right)^{2}  \varepsilon  \left( v(z,\theta, \omega) \right)^{2}   + \mathcal{O}\left( \varepsilon^{\min(1,2-2h)} \right).
\end{equation}

%\textcolor{red}{We emphasize that the above scale $k_1=\overline{k}_1\varepsilon^2$, i.e., $k_1=\mathcal{O}(\varepsilon^2)$, is crucial to our reconstruction scheme. The reason is that, for the eigenvalues $\lambda_n=\mathcal{O}(\varepsilon^2)$, the resonant regime is connected to the almost-vanishing nature of $1-\frac{\rho_1\omega^2\lambda_{n_0}}{k_1}$. Then, by assuming that $\rho_1$ and $\omega^{2}$ are both moderate constants, the only scale of $k_1$ that ensures the almost vanishing character of $1-\frac{\rho_{1}\omega^{2}\lambda_{n_0}}{k_1}$ is given by \eqref{Ctek1}, see \cite{Als} and section 4 in \cite{DGS} for more details. Moreover, the resonance regime will not occur for $\alpha \neq 2$, that is, the enhancement of the acoustic field is not justified. Thus we must take $\alpha=2$.}

The following remark is intended to provide the necessary clarification regarding the criteria governing the selection of the observation direction $\hat{x}$ in equation $(\ref{DiffUinfvinf})$, as this choice is not arbitrary but rather dictated by the underlying physical considerations.
\begin{remark}\label{COEF}
For other directions instead of the back-scattered direction in \eqref{DiffUinfvinf}, the information that we can recover is the product $v(z,-\hat{x},\omega)v(z,\theta, \omega)$ from which the extraction of $v(z,\theta, \omega)$ is tough. Hence \eqref{Heq} cannot be used to recover $k_0(\cdot)$ as the data $v(\cdot,\theta, \omega)$ is not available. So, as suggested in \cite{DGS}, we select the incident direction to be the back-scattered direction, and then \eqref{DiffUinfvinf} allows us to recover $\left(v(\cdot,\theta, \omega)\right)^2$, and then $v(\cdot,\theta, \omega)$ up to a sign.
\end{remark}

Obviously, in \eqref{BSD}, the field $v^{\infty}(-\theta, \theta, \omega)$ is independent of $z$, representing the far-field of the inhomogeneous background medium before injecting the droplet $D_z$. However, $u_z^{\infty}(-  \theta, \theta, \omega)$ from the droplet $D_z$ depends on $z$.  The identity  $(\ref{BSD})$  gives us an approximate scheme to compute $\left(v(z,\theta,\omega)\right)^{2}$ from the contrast
$u^{\infty}_z(-  \theta, \theta, \omega)-v^{\infty}(-  \theta, \theta, \omega)$, which will be taken as our inversion input for recovering $k_0(\cdot)$ in $\Omega$.

For given $n_{0} \in \mathbb{N}$, we will choose the incident frequency $\omega$ in the sequel such that
\begin{equation}\label{eq1039}
    \omega^{2}  =  \omega_{n_{0}}^{2}+  \mathcal{O}(\varepsilon^{h}) \quad \text{ for } 0 < h < 1,
\end{equation}
with $\omega_{n_{0}}^{2}$ being the eigen-frequency defined by \eqref{Def-Eigen-Freq}, where we take $\overline{k}_{1} = 1$, i.e., $\omega_{n_{0}}^{2}$ is given by
\begin{equation}\label{Eigfreq}
    \omega_{n_0}^2  =  \lambda_{n_0}^{-1}.
\end{equation}
Notice, in \eqref{eq1039} and \eqref{Eigfreq}, the small radius $0<\varepsilon \ll 1 $ and $n_0$ should ensure $\lambda_{n_0}^{-1}\gg \varepsilon^h$ for numerical computations. This condition is always satisfied as $\lambda_{n_0}  \sim 1$,  $h > 0$ and $\varepsilon$ is small. More precisely, using \eqref{eq1039} and \eqref{Eigfreq}, we are led to rewrite \eqref{BSD} as
\begin{equation}\label{CO1}
u^{\infty}_z(- \theta, \theta, \omega_{n_{0}})  - v^{\infty}(-  \theta, \theta, \omega_{n_0}) =   -  \frac{1}{4  \pi} \lambda_{n_{0}}^{-2}  \left( \int_B e_{n_{0}}(x)  dx \right)^{2}  \varepsilon^{1-h}  \left( v(z,\theta, \omega_{n_0}) \right)^{2}  +  \mathcal{O}\left( \varepsilon^{\min(1,2-2h)} \right).
\end{equation}

Next, to ensure that the first term on the right hand side of the above equation dominates the error term, i.e., $\varepsilon^{1-h} \gg  \varepsilon^{\min(1,2-2h)}$, we need to choose $h  \in  (1/2,1)$. Then, under the condition $\dfrac{1}{2} < h < 1$, the equation \eqref{CO1} becomes
\begin{equation}\label{CO2}
u^{\infty}_z(- \theta, \theta, \omega_{n_{0}})  - v^{\infty}(-  \theta, \theta, \omega_{n_0}) =   -  \frac{1}{4  \pi} \lambda_{n_{0}}^{-2}  \left( \int_B e_{n_{0}}(x)  dx \right)^{2}  \varepsilon^{1-h}  \left( v(z,\theta, \omega_{n_0}) \right)^{2}  +  \mathcal{O}\left( \varepsilon^{2-2h} \right).
\end{equation}
Besides, as measurements are rarely perfect, they are inevitably contaminated by noise fluctuations that obscure the true signal, e.g. the residual part $\mathcal{O}\left(\varepsilon^{2-2h}\right)$ in our formula \eqref{CO2}. Since the mentioned noise is typically small in magnitude compared to the underlying phenomenon, we can treat the measurement as the sum of a dominant, meaningful component and a residual, noisy one. By focusing on the dominant part and discarding the residual, we isolate the reliable information while filtering out the unreproducible irregularities that would otherwise distort our analysis. Hence, from \eqref{CO2}, we derive the following approximation
\begin{equation}\label{CO3}
u^{\infty}_z(-  \theta, \theta, \omega_{n_{0}})  - v^{\infty}(-  \theta, \theta, \omega_{n_0})  \approx  -  \frac{1}{4  \pi} \lambda_{n_{0}}^{-2}  \left( \int_B e_{n_{0}}(x)  dx \right)^{2}  \varepsilon^{1-h}  \left( v(z,\theta, \omega_{n_0}) \right)^{2}.
\end{equation}

Furthermore, when reconstructing an unknown function from measurements, our goal is not to recover perfect details but to capture the overall behavior. Even if our data are only an approximation, see \eqref{CO3}, it can be treated as exact within the context of the reconstruction process. Since the function itself is unknown, the sought of an exact representation is both impossible and unnecessary. Instead, we embrace the approximation as our working reality. This pragmatic approach allows us to build a stable and useful model from the information we have, without being paralyzed by its imperfections. Then the equation \eqref{CO3} can be rewritten as
\begin{equation}\label{BSD-1}
u^{\infty}_z(-  \theta, \theta, \omega_{n_{0}})  - v^{\infty}(-  \theta, \theta, \omega_{n_0})  =  -  \frac{1}{4  \pi} \lambda_{n_{0}}^{-2}  \left( \int_B e_{n_{0}}(x)  dx \right)^{2}  \varepsilon^{1-h}  \left( v(z,\theta, \omega_{n_0}) \right)^{2}.
\end{equation}
In the right-hand side of \eqref{BSD-1}, the determination of two constants $\lambda_{n_0}$ and $\int_B e_{n_{0}}(x) dx$ is of great importance. Based on
\cite[Theorem 4.2]{kalmenov2011boundary},  these two constants can be determined from the following result.

\begin{lemma}
For the Newtonian operator $N_B(\cdot)$ defined by \eqref{Def-N_B-Equa0708}, it follows that
\begin{equation}\label{EigValFractionalOrder}
\lambda_{n} = \mu_{n}^{-2}, \quad  n =1,\cdots,
\end{equation}
where $\mu_{n}$ are the positive roots of the transcendental equation
\begin{equation}\label{muFractionalOrder}
\textbf{J}_{\frac{1}{2}}
\left( \mu_{n} \right) + \frac{\mu_{n}}{2} \; \left(
\textbf{J}_{ - \frac{1}{2}}
\left( \mu_{n} \right)   -
\textbf{J}_{\frac{3}{2}}
\left( \mu_{n} \right) \right)  =  0,
\end{equation}
with
$\textbf{J}_{\nu}\left( \cdot \right)$  being the Bessel function of fractional order for $\nu \in \mathbb{R}$. The eigenfunctions corresponding to each eigenvalue $\lambda_{n}$ satisfy
\begin{equation}\label{EigFctFO}
\int_{B} e_{n}(x)  dx  =  4  \pi  \mu_{n}^{-\frac{5}{2}}  \int_{0}^{\mu_{n}} r^{\frac{3}{2}}
\textbf{J}_{\frac{1}{2}}
\left( r \right)  dr, \quad n=1,\cdots
\end{equation}
\end{lemma}

Thanks to this lemma, and by using the explicit expressions for Bessel functions $\textbf{J}_{-\frac{1}{2}}(x), \textbf{J}_{\frac{1}{2}}(x),$ and $ \textbf{J}_{\frac{3}{2}}(x)$,
the equation \eqref{muFractionalOrder} can be rewritten as
\begin{equation}\label{Eq112}
    \sin\left( \mu_{n} \right)  +  2   \mu_{n} \cos\left( \mu_{n} \right) =  0.
\end{equation}
Besides, the formula $(\ref{EigFctFO})$ can be simplified to
\begin{equation}\label{EigFctFON}
\int_{B} e_{n}(x)  dx   =  4  \sqrt{2 \pi}  \mu_{n}^{-\frac{5}{2}}  \int_{0}^{\mu_{n}} r  \sin(r)  dr\, \overset{(\ref{Eq112})}{=} 6   \sqrt{2 \pi}   \mu_{n}^{-\frac{5}{2}} \sin\left( \mu_{n}\right), \quad n=1,\cdots
\end{equation}
Then, for $\mu_{n}$ determined by $(\ref{Eq112})$, the eigenvalue $\lambda_{n}$ and $\int_{B} e_{n}(x)  dx$ can be computed from $(\ref{EigValFractionalOrder})$ and  \eqref{EigFctFON}, respectively.
The distribution of the first five roots of \eqref{Eq112} and the quantitative information about $\{\left(\mu_n,\lambda_n,\omega_n,\int_{B} e_{n}(x)dx \right): n=1,\cdots,5\}$ are shown in  Table 1.
%%%%%%%%%%%%%%%%%%%%%%%%%%%%%%
\begin{center}
\begin{table}[H]
\begin{tabular}{ |c|c|c|c|c|c| }
 \hline
 $n$ & 1 & 2 & 3 & 4 & 5  \\
 \hline
  $\mu_{n}$ & 1.8366 & 4.8158 & 7.9171 & 11.0408 & 14.1724 \\
 \hline
 $\lambda_{n}$ & 0.2965 & 0.0431 & 0.0160 & 0.0082 & 0.0050  \\
 \hline
  $\omega_{n}^{2}$ & 3.3726 & 23.2018 & 62.5 & 121.9512 &  200 \\
 \hline
  $\int_{B} e_{n}(x) \, dx $ & 3.1745 & -0.2939 & 0.0851 & -0.0371 & 0.0199  \\
  \hline
\end{tabular}
\caption{The distribution of  $\left\{\mu_{n}, \lambda_{n}, \omega_{n}^{2}, \int_{B} e_{n}(x)  dx \right\}_{n = 1}^{5}$.}
\label{Table1}
\end{table}
\end{center}
%%%%%%%%%%%%%%%%%%%%%%%%%%%%%%%%%%%%%
Based on Table \ref{Table1}, we can estimate the following constant
\begin{equation}\label{DefCn0}
    C_{n_{0}}  :=  \frac{1}{4  \pi}   \lambda_{n_{0}}^{-2} \left( \int_{B} e_{n_{0}}(x) dx \right)^{2}  \varepsilon^{1-h},
\end{equation}
which, combined with $(\ref{BSD-1})$, gives us
\begin{equation}\label{BSD-1-Algo}
 v^{\infty}(-  \theta, \theta, \omega_{n_0})  -  u^{\infty}_z(-  \theta, \theta, \omega_{n_{0}})  = C_{n_{0}} \,   \left( v\left(z,\theta, \omega_{n_0}\right) \right)^{2}.
\end{equation}

\medskip
We numerically recover the bulk modulus function $k_{0}(z)$ within the bounded region $\Omega$ using the data $v^{\infty}(-  \theta, \theta, \omega_{n_0})  -  u^{\infty}_{z}(-  \theta, \theta, \omega_{n_{0}})$. It should be noted that, due to \eqref{BSD-1-Algo}, this scheme is only applicable to points $z_0 \in \Omega$ satisfying $v(z_0,\theta,\omega_{n_0}) \neq 0$. Indeed, if $v(z_0,\theta,\omega_{n_0}) = 0$, the quantity $\dfrac{1}{k_0(z_0)}$ will physically disappear  in \eqref{Heq}.
However, the analyticity of $v(\cdot)$ in $\Omega$ ensures that the measure of zero points of $v(\cdot,\theta,\omega_{n_0})$ is zero, so $k_0(z_0)$ for $z_0$ being the zero points of $v(\cdot,\theta,\omega_{n_0})$ can be obtained by the limitation process $z^*\to z_0$ for $z^*\in\Omega$ satisfying $v(z^*,\theta,\omega_{n_0})\not=0$ with known recovered $k_0(z^*)$. Due to this reason, we assume that $v(\cdot,\theta,\omega_{n_0})\not=0$ in the whole domain $\Omega$ for the simplicity of statement.

Building upon the relations established above, the algorithm for the numerical reconstruction of the bulk modulus function $k_0(\cdot)$ in $\Omega$, as introduced in \cite{DGS}, is given by:
\begin{algo}\label{AlgoScheme}
\begin{enumerate}
\item[]
    \item \label{Stepa} \textbf{Step (1)}.  We start by collecting our inversion input
    \begin{equation}\label{KK}
       \xi(z) :=  \xi(z,\theta, \omega_{n_{0}})  :=   v^{\infty}(-  \theta, \theta, \omega_{n_0}) -  u^{\infty}_z(-  \theta, \theta, \omega_{n_{0}}),
    \end{equation}
    with $z\in\Omega$, which corresponds to the left hand side of \eqref{BSD-1-Algo}.

    \item \label{Stepb} \textbf{Step (2)}. Based on \eqref{DefCn0}, \eqref{BSD-1-Algo},  and \eqref{KK}, we compute $v(z,\theta,\omega_{n_{0}})$, which is the point-wise value of the total field  before injecting the droplet $D_{z}$ into $\Omega$.

     \item \label{Stepc} \textbf{Step (3)}. We use numerical differentiation techniques (mollification method) to compute $\Delta v(z,\theta,\omega_{n_{0}})$ from the already reconstructed value $v(z,\theta,\omega_{n_{0}})$ in $\Omega$ ,as shown in \textbf{Step $(\ref{Stepb})$}.

     \item \label{Stepd} \textbf{Step (4)}. We use the Helmholtz equation \eqref{Heq} to obtain a set of point-wise values of $k_{0}(\cdot)$  by
     \begin{equation}\label{liu03-01}
        \frac{1}{k_0(z)}  =   - \frac{1}{\omega_{n_0}^2}  \frac{\Delta v(z,\theta,\omega_{n_0})}{v(z,\theta,\omega_{n_0})},\quad z \in \Omega,
     \end{equation}
    where, in the right-hand side, $v(z,\theta,\omega_{n_0})$ is already reconstructed from \textbf{Step (\ref{Stepb})}, and $\Delta v(z,\theta,\omega_{n_0})$ is already reconstructed from \textbf{Step (\ref{Stepc})}. Hence, by moving the droplet $D_{z}$ in $\Omega$, we can recover $\left\{k_{0}(z_{j})\right\}_{j=1}^{N}$, the point-wise values of $k_{0}(\cdot)$.
\end{enumerate}
\end{algo}

Furthermore, based on the \textbf{Algorithm \ref{AlgoScheme}}, the continuous values of $k_0(\cdot)$ can be finally approximated from  $\left\{k_{0}(z_{j})\right\}_{j=1}^{N}$ obtained in \textbf{Step (\ref{Stepd})} from the expansion
\begin{equation}\label{coro-Eq2.8}
    k_{0}(x)  =  \sum_{j=1}^{N} \gamma_{j}  f_{j}(x), \quad x\in \Omega,
\end{equation}
where $(f_{1}(\cdot), \cdots, f_{N}(\cdot))$ is the specified radial basis functions, for example, see \textbf{Subsection \ref{SubSectionUnperturbed}} for the choice of the radial basis functions. The expansion coefficient $\vec{\gamma}:=(\gamma_{1}, \cdots, \gamma_{N})^T \in  \mathbb{C}^{N}$ can be determined from the following linear algebraic system
\begin{equation*}\label{coro-LZEq1037}
(f_i(z_j))_{N\times N}\cdot\vec{\gamma}=\vec{k_0}
\end{equation*}
with already reconstructed vector $\vec{k_0}:=(k_0(z_1),\cdots,k_0(z_N))^T\in \mathbb{C}^N$, see \textbf{Step (\ref{Stepd})}.

Some necessary clarifications to the implementation of   \textbf{Algorithm \ref{AlgoScheme}} are stated as follows.

\begin{itemize}

\item Concerning the selected index.

Since the index $n_{0} \in \mathbb{N}$ in \eqref{KK} is arbitrary,  we take $n_{0} = 1$  in \eqref{KK} to obtain
\begin{equation}\label{BSDI}
\xi(z)  =  \xi(z,\theta, \omega_{1})  =  v^{\infty}(- \theta, \theta, \omega_1)  -  u_z^{\infty}(- \theta, \theta, \omega_1)   \overset{(\ref{BSD-1-Algo})}{=}  C_{1} \left( v(z,\theta, \omega_1) \right)^{2}, \quad z \in \Omega,
\end{equation}
without loss of generality, where $\omega_1   \overset{(\ref{Eigfreq})}{=}  \lambda_1^{-1/2}  \overset{( \ref{EigValFractionalOrder})}{=}  \mu_{1}   \overset{\textbf{Table \ref{Table1}}}{\simeq}  1.8366$, and
\begin{equation}\label{KKWHW}
    C_{1}  \overset{(\ref{DefCn0})}{=}  \frac{1}{4  \pi} \lambda_{1}^{-2}  \left(\int_{B} e_{1}(x)  dx \right)^{2}  \varepsilon^{1-h} \overset{\textbf{Table \ref{Table1}}}{\simeq} \frac{1}{4  \pi}  \left( 0.2965 \right)^{-2}  \left( 3.1745\right)^{2}  \varepsilon^{1-h}
\end{equation}
with $h >1/2$. Under this configuration, the equation \eqref{liu03-01} becomes
\begin{equation}\label{liu03}
\frac{1}{k_0(z)}  = - \frac{1}{\omega_{1}^2} \frac{\Delta v(z,\theta,\omega_{1})}{v(z,\theta,\omega_{1})},\quad z\in \Omega.
\end{equation}

\item  Concerning the numerical differentiation involved in \textbf{Step (\ref{Stepc})}.

To handle the numerical differentiation of $v(z,\theta,\omega_{1})$ such as $\Delta v(z,\theta,\omega_{1})$, see for instance \textbf{Step $(\ref{Stepc})$}, we express the complex-valued inversion input data as
\begin{equation}\label{liu0201}
\xi(z)   = C_{1} \left(v(z, \theta, \omega_{1})\right)^2, \quad z \in \Omega,
\end{equation}
see $(\ref{BSDI})$, with $C_{1}$ is estimated by $(\ref{KKWHW})$. Then, from $(\ref{liu0201})$, we deduce that
\begin{equation}\label{Equa1129}
    v(z, \theta, \omega_{1})  =  \frac{\sqrt{\xi(z)}}{\sqrt{C_{1}}}, \quad z \in \Omega.
\end{equation}
Hence, for $z \in \Omega$, from \eqref{liu03} and \eqref{Equa1129}, it follows that
\begin{equation}\label{rek0}
\frac{1}{k_0(z)} =  - \frac{1}{\omega_{1}^2}  \frac{\Delta \sqrt{\xi(z)}}{\sqrt{\xi(z)}}=-  \frac{1}{\omega_{1}^2}  \left( \frac{\Delta\xi(z)}{2\xi(z)}-\frac{|\nabla\xi(z)|^2}{4\xi(z)^2}\right).
\end{equation}

\end{itemize}

A very important observation in \eqref{rek0} is  that, once we determine $\xi(z)$ from  \eqref{BSDI}, the right-hand side of  \eqref{rek0} is independent of $C_1$, i.e.,  we do not need $C_1$ approximated by \eqref{KKWHW} for our reconstruction. On the other hand, the right-hand side of  \eqref{rek0} requires numerical differentiations, which are  ill-posed for $\xi(z)$ given by \eqref{BSDI}, since
$u_z^{\infty}\left(-\theta, \theta, \omega_{1}\right)$ are approximated by solving truncated integral equations for small $0 <\varepsilon \ll 1$, and consequently $\xi(z)$ is of unavoidable noise.  To address this issue, we utilize the mollification method to deal with the ill-posedness, see \cite{mol1,Mur}
for related regularization schemes. More precisely, for given  $\xi(\cdot)$ on a domain $\mathfrak{U}$ containing the point $z$ where we want to compute its Laplacian, we compute its second order partial derivatives along three directions, respectively. We firstly  take $\mathfrak{U}\subset \mathbb{R}^3$ to be a cube and  introduce
\begin{equation*}
  \eta(x):= \begin{cases}
			\textbf{c} \exp\left(\dfrac{1}{\left| x\right|^2-1}\right), & \text{if $\left\vert x \right\vert < 1$,}\\
            0, & \text{if $\left\vert x \right\vert \geq 1$,}
		 \end{cases}
\end{equation*}
with
     $\textbf{c}:= \left( \int^{1}_{-1}e^{\frac{1}{\left|x \right|^2-1}}dx \right)^{-1}$, which is a $C^{\infty}_0(\mathbb{R}^1)$ function satisfying $\int_{\mathbb{R}^1}\eta(x)dx=1$. For kernel function $\eta_{\delta}(x):=\frac{1}{\delta}\eta(\frac{x}{\delta})$ with specified $\delta>0$, if $f\;:U \subset \mathbb{R}\to \mathbb{R}$ is locally integrable, define its mollification \cite{EvansBook}
        \begin{equation*}
            f^{\delta}(x) := (\eta_{\delta}*f)(x) :=\int_{U}\eta_{\delta}(x-y)f(y)dy=\int_{-\delta}^{\delta}\eta_{\delta}(y)f(x-y)dy, \qquad x \in   U_{\delta},
        \end{equation*}
        where $U_{\delta}:= \left\{x \in U \mid  \hbox{dist}(x,\partial U)>\delta \right \}$. Additionally, from \cite[Chapter 2, Page 37]{Adam}, we have
     \begin{equation}\label{FARNJB}
         D^{\alpha}f^{\delta}(x) =  \int_{U} D^{\alpha}\eta_{\delta}(x-y)f(y)  dy, \quad x \in U_{\delta}\subset U\subset \mathbb{R},
     \end{equation}
     where $D=\frac{d}{dx}$, is the differentiation operator,  and $\alpha=1,2$ is the order of derivative. According to \cite[Theorem 4.1.1]{molill}, the formula \eqref{FARNJB} makes sense under the condition $\left\Vert D^{\alpha}f \right\Vert_{\mathbb{L}^{p}\left( U_{\delta}\right)} \le \mathrm{M}(\alpha, p, \delta)$, for $1 \le p \le \infty$, with some  positive constant $\mathrm{M}(\alpha, p, \delta)$, and it holds that
        \begin{equation}\label{VUPR}
           \lim_{\delta \to  0} \left\Vert D f^{\delta} - Df \right\Vert_{\mathbb{L}^{p}\left( U_{\delta} \right)} \to 0 \quad \text{and} \quad  \lim_{\delta  \to  0} \left\Vert D^2f^{\delta} - D^2f \right\Vert_{\mathbb{L}^{p}\left( U_{\delta} \right)} \to 0.
        \end{equation}

The following remark is necessary to clarify the $\mathbb{L}^{p}\left( U_{\delta} \right)$-convergence properties of $(\ref{VUPR})$.

\begin{remark}
The convergence property \eqref{VUPR} is stated for exact data $f$.
In practice we always have noisy measurements $f^\tau$ satisfying
$$\|f^\tau-f\|_{L^p(U_{\delta})}\le \tau$$
with noise level $\tau>0$. In this situation, we have to take the mollification process  $f^{\delta,\tau}:=\eta_\delta*f^\tau$ for noisy data $f^\tau$. Then, for any index $\alpha=1,2$,
we make the decomposition
\[
D^\alpha f^{\delta,\tau} - D^\alpha f
=
D^\alpha(\eta_\delta*(f^\tau-f))
+
\big(D^\alpha f^\delta-D^\alpha f\big).
\]
The first term $D^\alpha(\eta_\delta*(f^\tau-f))$ is called the noise propagation term, and the second term $D^\alpha f^\delta-D^\alpha f$ is the mollification
bias term, which tends to $0$ as the mollification parameter $\delta\to 0$ from \eqref{VUPR}.
By the commutation of differentiation and convolution and Minkowski's inequality\cite[Theorem 1.2.10]{Minkowski}, we have
\[
\|D^\alpha(\eta_\delta*(f^\tau-f))\|_{L^p(U_\delta)}
=
\|(D^\alpha\eta_\delta)*(f^\tau-f)\|_{L^p(U_\delta)}
\le \|D^\alpha\eta_\delta\|_{L^1(U_\delta)}\,\|f^\tau-f\|_{L^p(U_\delta)}
\le C\tau\delta^{-\alpha},
\]
where $C=\|D^\alpha\eta\|_{L^1(U_{\delta})}$ from $\|D^\alpha\eta_\delta\|_{L^1(U_{\delta})}=\delta^{-\alpha}\|D^\alpha\eta\|_{L^1(U_{\delta})}$.
So a possible strategy for choosing the mollification parameter  $\delta=\delta(\tau)$ in terms of noise level $\tau$ is
\begin{equation}
\delta,\; \tau\delta^{-\alpha}\to 0\quad \hbox{as }\tau\to 0
\end{equation}
from which $D^\alpha f^{\delta(\tau),\tau}\approx D^\alpha f$ holds as $\tau\to 0$.
\end{remark}

To compute $\nabla\xi(z)$ (respectively. $\Delta\xi(z)$), we start by using $(\ref{FARNJB})$ to compute all the partial derivatives related to $\nabla\xi^{\delta}(z)$ (respectively. $\Delta\xi^{\delta}(z)$) to form
\begin{eqnarray*}
    \nabla\xi^{\delta}(z)=\left(\frac{\partial\xi^{\delta}}{\partial z_1}(z), \frac{\partial\xi^{\delta}}{\partial z_2}(z), \frac{\partial\xi^{\delta}}{\partial z_3}(z)\right) \quad \left(\text{respectively.} \; \Delta\xi^{\delta}(z) = \sum_{\ell=1}^{3} \frac{\partial^2\xi^{\delta}}{\partial z_\ell^2}(z) \right)
\end{eqnarray*}
with specified mollification constant $\delta>0$, and then thanks to $(\ref{VUPR})$ we deduce an approximations of $\nabla\xi(z)$ (respectively. $\Delta\xi(z)$). To do this, for fixed $(z_2,z_3)$, we assume that $\xi(\cdot,z_2,z_3): U \to \mathbb{R}$ to be a locally integrable function, then we have
\begin{eqnarray*}
        \xi^{\delta }(z_1,z_2,z_3):= \left( \eta_{\delta}*\xi(\cdot,z_2,z_3)\right)(z_1)=\int_{U}\eta_{\delta}(z_1-t)\xi(t,z_2,z_3)dt, \quad z_1\in U_{\delta}.
\end{eqnarray*}
Thus, from $(\ref{FARNJB})$, it holds that
\begin{equation*}
    \partial_{z_1}\left(\xi^{\delta }(z_1,z_2,z_3)\right) = \int_{U}\partial_{z_1}\eta_{\delta}(z_1-t)\xi(t,z_2,z_3)dt, \quad z_1\in U_{\delta},
\end{equation*}
and, in a similar manner,
\begin{eqnarray*}
        \partial_{z_2}\left(\xi^{\delta }(z_1,z_2,z_3)\right)&=&\int_{U}\partial_{z_2}\eta_{\delta}(z_2-t)\xi(z_1,t,z_3)dt, \quad z_2\in U_{\delta},\\
        \partial_{z_3}\left(\xi^{\delta }(z_1,z_2,z_3)\right)&=&\int_{U}\partial_{z_3}\eta_{\delta}(z_3-t)\xi(z_1,z_2,t)dt, \quad z_3\in U_{\delta}.
\end{eqnarray*}
Hence, as explained above, we can compute an approximation to $\nabla\xi(z)$ by computing
\begin{equation}\label{grexi}
     \nabla\xi^{\delta}(z) \, = \, \left(\frac{\partial\xi^{\delta}}{\partial z_1}, \frac{\partial\xi^{\delta}}{\partial z_2}, \frac{\partial\xi^{\delta}}{\partial z_3}\right)(z).
\end{equation}
Similarly, an approximation to $\Delta \xi(z)$ by
\begin{equation}\label{lapxi}
    \Delta\xi^{\delta}(z) =\frac{\partial^2\xi^{\delta}}{\partial z_1^2}(z) + \frac{\partial^2\xi^{\delta}}{\partial z_2^2}(z)+\frac{\partial^2\xi^{\delta}}{\partial z_3^2}(z),
\end{equation}
can be derived by computing
\begin{equation*}
        \partial_{z_1}^2\left(\xi^{\delta }(z_1,z_2,z_3)\right)=\int_{U}\partial_{z_1}^2\eta_{\delta}(z_1-t)\xi(t,z_2,z_3)dt, \quad z_1\in U_{\delta},
\end{equation*}
and
\begin{eqnarray*}
 \partial_{z_2}^2\left(\xi^{\delta }(z_1,z_2,z_3)\right)   &=& \int_{U} \partial_{z_2}^2\eta_{\delta}(z_2-t)\xi(z_{1},t,z_3)\, \, dt, \quad z_2\in U_{\delta}, \\
  \partial_{z_2}^2\left(\xi^{\delta }(z_1,z_2,z_3)\right)  &=& \int_{U} \partial_{z_3}^2\eta_{\delta}(z_3-t) \xi(z_{1},z_2,t)\, \, dt, \quad z_3\in U_{\delta}.
\end{eqnarray*}

Before proceeding to the reconstruction of the pointwise value of bulk modulus function $k_{0}(z)$, we note that the above technique applies the one-dimensional mollifier $\eta_{\delta}(\cdot)$ component-wise to smooth the 3D inversion input data $\xi(z)$ in preparation for the 3D Laplacian operation, which is an approach readily implemented numerically within a cube. However, this particular regularization scheme may introduce an undesired directional bias, warranting further investigation into a general-form smoothing mollifier for the 3D Laplacian to mitigate such artifacts. Addressing this limitation through a general-form smoothing mollifier for the 3D Laplacian remains an important direction for future investigation.
%Here we apply  the 1Dmollifier $\eta_{\delta}$  to smoothen the 3D inversion input data $\xi(z)$ component-wisely for 3D Laplacian operation, which is easily to be implemented numerically in a cube. However, such a special regularizing scheme may introduce extra directional bias. To improve such an extra error, the smoothing mollifier in general form for 3D Laplacian operation should be further studied.}

Thus, by returning back to $(\ref{rek0})$, we deduce
\begin{equation}\label{XC-Equa0850}
\frac{1}{k_0(z)}  = -  \frac{1}{\omega_{1}^2}  \left( \frac{\Delta\xi^{\delta}(z)}{2\xi(z)}-\frac{|\nabla\xi^{\delta}(z)|^2}{4\xi(z)^2}\right) - \frac{1}{\omega^{2}}  Error(z), \quad z \in \mathfrak{U}_{\delta},
\end{equation}
where the error term $Error(z)$ is given by
\begin{eqnarray*}
 Error(z) &:=&  \frac{\Delta\left(\xi-\xi^{\delta}\right)(z)}{2\xi(z)}  -   \frac{|\nabla\left(\xi -\xi^{\delta}\right)(z)|^2}{4 \xi(z)^2}- \\ & & \frac{\langle Re\left(\nabla \xi (z) \right);Re\left(\nabla \xi (z)-\nabla \xi^{\delta}(z) \right)\rangle  +  \langle Im\left(\nabla \xi (z) \right);Im\left(\nabla \xi (z)-\nabla \xi^{\delta}(z) \right)\rangle}{2  \xi(z)^2},
\end{eqnarray*}
which, by using $(\ref{VUPR})$, can be estimated as
\begin{equation*}
\left\Vert Error \right\Vert_{\mathbb{L}^{p}\left( \mathfrak{U}_{\delta} \right)} \lesssim \left\Vert \Delta \xi^{\delta} - \Delta \xi \right\Vert_{\mathbb{L}^{p}\left( \mathfrak{U}_{\delta} \right)} + \left\Vert \nabla \xi^{\delta} - \nabla \xi \right\Vert_{\mathbb{L}^{2p}\left( \mathfrak{U}_{\delta} \right)}  \to 0 \;\; \text{as} \;\; \delta \to 0.
\end{equation*}
Hence, by referring to \eqref{XC-Equa0850}, we have
\begin{equation}\label{RF}
 \dfrac{1}{k_{0}(z)}\approx -  \frac{1}{\omega_{1}^2} \left( \frac{\Delta\xi^{\delta}(z)}{2\xi(z)}-\frac{|\nabla\xi^{\delta}(z)|^2}{4\xi(z)^2}\right), \quad z \in \mathfrak{U}_{\delta}\subset \mathfrak{U}\subset \mathbb{R}^3,
\end{equation}
where $\xi(z)$ is given by \eqref{BSDI}, while $\nabla\xi^{\delta}(z)$, and $\Delta\xi^{\delta}(z)$ are given by \eqref{grexi}, and \eqref{lapxi}, respectively.
\medskip
\newline
We conclude this section by the following explanations.
\begin{remark}
The following points are in order.
\begin{enumerate}

    \item Due to the quotient structure in the right-hand side of \eqref{liu03}, the constant difference of the argument when we determine the complex-valued $v(z, \theta, \omega_{1})$ from $\left( v(z, \theta, \omega_{1}) \right)^2$ has nothing to do with the reconstruction of $k_0(z)$.

    \item In \textbf{Step $(\ref{Stepc})$}, we use the mollification method to calculate $\Delta v(z, \theta, \omega_{1})$ from $v(z, \theta, \omega_{1})$, but there are other methods that can also be used, such as finite difference methods, the Tikhonov regularization, the level set methods, and the conjugate gradient techniques, as shown in \cite{wave,con, FDM1,FDM2}, and the references therein.

    \item The point-wise values of the bulk modulus function $k_{0}(\cdot)$ can only be recovered in regions where the total field $v(\cdot, \theta, \omega_{1})$ is not constant, see \eqref{Heq} and \eqref{liu03-01}.

    \item  To recover $k_{0}(\cdot)$ from its discrete values, we have opted for interpolation with radial basis functions in \eqref{coro-Eq2.8}. Obviously, there are other methods that can be employed such as the spline interpolation.

    \item We apply one single droplet for our reconstruction scheme. From the practical situations,  it
      is more reasonable to do imaging using multiple droplets. In principle, the reconstruction algorithm established in this work is also applicable to the multiple droplets case. For a more detailed discussion in the context of the electromagnetic model, see \cite{Cao}.

\end{enumerate}
\end{remark}
%%%%%%%%%%%%%%%%%%%%%%%%%%%%%%%%%%%%%%%%%%%%%%%%%%%%%%%
%%%%%%%%%%%%%%%%%%%%%%%%%%%%%%%%%%%%%%%%%%%%%%%%%%%%%%
\section{The direct scattering for inhomogeneous medium}\label{SectionII}

The goal of this section is to derive the far-field patterns of the scattered waves for inhomogeneous medium with and without the droplet denoted by $D_{z}$ with small size and high contrasting bulk modulus, see \eqref{DzepsB}, \eqref{defrhodefk} and \eqref{Ctek1}, which is necessary for our inversion scheme, see \textbf{Step $(\ref{Stepa})$}.
%The direct scattering of acoustic wave by inhomogeneous medium with droplets itself is of great importance, since the droplet $D_{z}$ with small size and high contrasting bulk modulus, embedded in inhomogeneous medium, will change the scattering behaviors compared with general inhomogeneous medium.
%On the other hand, the solution to the direct scattering process will yield simulation of the far-field patterns which are taken as inversion input for our reconstruction scheme.

Our basic scheme, to determine the far-field waves $u_z^{\infty}(\cdot,\cdot,\cdot)$ and $ v^{\infty}(\cdot,\cdot,\cdot)$, is to firstly generate the total waves $u(\cdot,\cdot,\cdot)$ and $v(\cdot,\cdot,\cdot)$ in  $B({\bf 0}, 1)$ for known $k_{0}(\cdot)$ by simulation process, and then to compute $u_z^{\infty}(\cdot,\cdot,\cdot)$ and $ v^{\infty}(\cdot,\cdot,\cdot)$ by their integral expressions, respectively. Consequently, we will fix the artificial domain $\Omega \equiv B({\bf 0},1)$ where we want to recover  $k_0(\cdot)$, without loss of generality. In fact, if we want to recover the inhomogeneous function $k_0(\cdot)$ in general domain $\Omega$ not necessarily in $B({\bf 0},1)$, we can always take $B({\bf 0}, R_0)$ with radius $R_0$ large enough such that $\Omega \subset B({\bf 0}, R_0)$. For this domain $B({\bf 0}, R_0)$, our implementation schemes for both direct and inverse scattering in $B({\bf 0}, 1)$ are still applicable, except for the introduction of the amplifying factor $R_0$ in the computational process.

\subsection{Generated  fields for the unperturbed medium}\label{SubSectionUnperturbed}

We recall from  \eqref{Heq}  that, the total field $v(\cdot,\theta,\omega)$ inside $\Omega  \equiv  B({\bf 0}, 1)$ for fixed $(\theta,\omega)$, before injecting the droplet $D_z$, is the solution to
\begin{equation}\label{HeqI}
    {\Delta} v(x,\theta, \omega)  +  \frac{\omega^{2}}{k_{0}(x)}  v(x,\theta, \omega)  =  0, \quad x \in \mathbb{R}^{3}, \; \theta \in \mathbb{S}^{2}.
\end{equation}
We decompose our simulation process into two steps. The first step is to compute the total field $v(\cdot,\theta,\omega)$, solution of  $(\ref{HeqI})$, while the second step is to compute its corresponding far-field $v^{\infty}(\cdot,\theta,\omega)$ in $\mathbb{S}^2$.

1. Computation of the total field $v(\cdot,\theta,\omega)$.

    To do this, using  \eqref{PDE-u-inc}, we can rewrite \eqref{HeqI} as
    \begin{equation*}
    {\Delta} v(x,\theta, \omega)  +  \frac{\omega^{2}}{k_{0}(x)}  v(x,\theta, \omega)  =  {\Delta} u^{inc}(x,\theta, \omega)  + \omega^{2} u^{inc}(x,\theta, \omega), \quad x \in \mathbb{R}^{3}, \; \theta \in \mathbb{S}^{2}.
\end{equation*}
    Multiplying the above equation by the fundamental solution
    \begin{equation}\label{FundSolHelmoltzEq}
        \Phi_{\omega}(x,y) :=   \dfrac{e^{i \omega \left\vert x - y \right\vert}}{4 \pi  \left\vert x - y \right\vert}, \quad x \neq y
    \end{equation}
to the Helmholtz equation, which meets  the equation  $\Delta \Phi_{\omega}(x,y)  +  \omega^{2}  \Phi_{\omega}(x,y)   =   -  \delta(x,y)$
    %\begin{equation}\label{FundSolHel}
    %    \Delta \Phi_{\omega}(x,y)  +  \omega^{2}  \Phi_{\omega}(x,y)   =   -  \delta(x,y)
    %\end{equation}
    in the distributional sense,  we derive the following Lippmann-Schwinger equation
    \begin{equation}\label{Eq0823}
        v(x, \theta, \omega)  -  \omega^{2}  \int_{B({\bf 0},1)} \left( \frac{1}{k_{0}(y)} -  1 \right)  \Phi_{\omega}(x,y)  v(y,\theta,\omega)  dy  =  u^{inc}(x,\theta, \omega), \quad x \in \mathbb{R}^{3},
    \end{equation}
in terms of the radiation conditions at infinity. From \eqref{Eq0823}, we observe that $v(\cdot,\theta,\omega)$ in $\mathbb{R}^{3}$ is completely determined by  $v(\cdot,\theta,\omega)$ in $B({\bf 0},1)$. Thus, the computation of $v(\cdot,\theta,\omega)$ in $B({\bf 0},1)$ is of great importance. To this end, by restricting the Lippmann-Schwinger equation \eqref{Eq0823} into $B({\bf 0},1)$, and denoting its corresponding solution by $v(\cdot, \theta, \omega)|_{B({\bf 0},1)}$, we get for $ x \in B({\bf 0},1)$ that
        \begin{equation}\label{liu04}
        v(x, \theta, \omega)|_{B(0,1)}  -  \omega^{2}  \int_{B({\bf 0},1)} \left( \frac{1}{k_{0}(y)} -  1 \right)  \Phi_{\omega}(x,y)  v(y, \theta, \omega)|_{B({\bf 0},1)}  dy = u^{inc}(x,\theta, \omega).
    \end{equation}
Once we solve \eqref{liu04}, which is a Fredholm integral equation of the second kind, the total field $v(\cdot, \theta, \omega)$ in the whole space $\mathbb{R}^3$ can be computed from \eqref{Eq0823}.
\medskip
\newline
There are many schemes in the literature to solve the Lippmann-Schwinger equation \eqref{liu04}, such as the Born and Rytov approximations  for linearized models \cite{Born, Rytov}, the contrast source-inversion method \cite{con_source}, and the recursive Born approximation for nonlinear model with heavy computational costs \cite{ReBorn}. The challenging numerical issue is the computation of the volume integral with weak singular kernel, due to the fundamental solution \eqref{FundSolHelmoltzEq} appearing on the left-hand side of \eqref{liu04}. To deal with this issue, according to \cite{DRM92,DRBEM}, the dual reciprocity method (DRM) is the most commonly used method for converting volume integrals into surface integrals based on the Green's formula. Then the boundary element method (BEM) can be applied to compute the corresponding surface integrals \cite{BEMBrebbia2, BEMDing}, which weakens the computational cost for the Lippmann-Schwinger equation \eqref{liu04}.
% BEMBrebbia1,
To apply the DRM to solve \eqref{liu04} efficiently, we firstly introduce the functions $v^{\star}(\cdot, \theta, \omega)$ and $u^{\star,inc}(\cdot, \theta, \omega)$ defined by
\begin{equation}\label{defvstar}
v^{\star}(x, \theta, \omega)  :=  \left( \frac{1}{k_{0}(x)} -  1 \right) v(x, \theta, \omega), \quad x\in B({\bf 0},1),
\end{equation}
\begin{equation}\label{Uincstar}
u^{\star,inc}(x, \theta, \omega) :=  \left( \frac{1}{k_{0}(x)} -  1 \right) u^{inc}(x, \theta, \omega), \quad x\in B({\bf 0},1),
\end{equation}
respectively. Then, multiplying both sides of \eqref{liu04} by  $\left( \dfrac{1}{k_{0}(\cdot)}  -  1 \right)$, using  \eqref{defvstar} and \eqref{Uincstar}, we obtain
 \begin{equation}\label{liu05}
         v^{\star}(x, \theta, \omega)  -  \omega^{2}
        \left( \frac{1}{k_{0}(x)}  -  1 \right) \int_{B({\bf 0},1)}   \Phi_{\omega}(x,y)   v^{\star}(y, \theta, \omega)  dy  =  u^{\star,inc}(x, \theta, \omega), \quad  x \in B({\bf 0},1).
    \end{equation}

Now, inspired by \cite[section 1]{Zhe-Wang}, the volume integral appearing on the left hand side of \eqref{liu05} can be equivalently transformed to a surface integral if and only if the density function over the volume integral, i.e., $v^{\star}(\cdot, \theta, \omega)$, is a source of a given Helmholtz problem. To handle this transformation, we will use the basic concept of DRM, which consists of approximating $v^{\star}(\cdot, \theta, \omega)$ by a linear combination of radial basis functions in the form
\begin{equation}\label{Vexpres}
    v^{\star}(x, \theta, \omega)  =  \sum_{k=1}^{n} \alpha_{k}  f_{k}(x),\quad x\in B({\bf 0},1),
\end{equation}
with $\left(\alpha_{1}, \cdots, \alpha_{n}\right)\in \mathbb{C}^n$ the expansion coefficients to be determined,  and $\left(f_{1}(\cdot), \cdots, f_{n}(\cdot)\right)$ is the suitably chosen radial basis functions,
see \cite{DRBEM, KR, GCBP}. For notations' simplicity, we omit the dependency of $\alpha _{k} = \alpha_{k}(\theta,\omega)$ on $(\theta,\omega)$. Clearly, if each basis function $f_{k}(\cdot)$ is chosen as a source term for a given Helmholtz problem, $v^{\star}(\cdot, \theta, \omega)$ can be approximated by a function which is a source for a given Helmholtz problem. There are already thorough studies on the choice strategy for the radial basis function in the DRM, see \cite{DRBEM, KR, GCBP} and the references therein, such as
$$f_{k}(\cdot)  =  1, \quad f_{k}(\cdot) =  1  +  \left\vert \cdot - x_{k} \right\vert,\quad f_{k}(\cdot)  = e^{\left\vert \cdot - x_{k} \right\vert^{2}},\quad f_{k}(\cdot)  =  \left\vert \cdot - x_{k} \right\vert^{2}\ln \left(\left\vert \cdot - x_{k} \right\vert \right),$$
where $\{x_{k}\in B({\bf 0},1): k=1,\cdots,n\}$ is a set of randomly collocated points. Now, by plugging  \eqref{Vexpres}  into  \eqref{liu05}, we obtain
\begin{equation}\label{Eq0936II}
        \sum_{k=1}^{n} \alpha_{k}  \left[  f_{k}(x)  -\omega^{2}  \left( \frac{1}{k_{0}(x)}-1 \right)\int_{B(0,1)}   \Phi_{\omega}(x,y) f_{k}(y)  dy  \right] = u^{\star,inc}(x, \theta, \omega), \quad x \in B({\bf 0},1).
    \end{equation}
\bigskip
To numerically handle the integral appearing on the left hand side of \eqref{Eq0936II} given by
\begin{equation}\label{Equa0500}
    u_{k}(x)  :=  \int_{B({\bf 0},1)}   \Phi_{\omega}(x,y) f_{k}(y)  dy, \quad x \in B({\bf 0},1),
\end{equation}
it is crucial to select a good radial basis function $f_{k}(\cdot)$ to obtain  good numerical quadrature  \eqref{Equa0500}. In many instances, as demonstrated in \cite{DRBEM}, the radial basis function
\begin{equation}\label{RBF}
     f_{k}(x) := 1 + \left\vert x  -  x_{k} \right\vert, \quad x \in B({\bf 0},1), \; 1 \leq k \leq n,
\end{equation}
where $x_k \in B({\bf 0},1)$, is the most frequently used one due to its satisfactory outcomes.
Thanks to \cite[Proposition 1.1]{Zhe-Wang}, and the expression of $u_{k}(\cdot)$, defined by \eqref{Equa0500}, we know that $f_{k}(\cdot)$ will be a source term for the following Helmholtz equation
\begin{equation}\label{PT-Equa0918}
    \Delta u_{k}(x)  + \omega^{2}  u_{k}(x) \, = \, - \, f_{k}(x), \quad x \in B({\bf 0},1),
\end{equation}
and $u_{k}(\cdot)$ satisfies the following boundary condition\footnote{$p.v$ stands for the Cauchy principal value.} for  $x \in \partial B({\bf 0},1)$ that
 \begin{equation*}\label{BCuk}
    \frac{1}{2}  u_{k}(x)  + \int_{\partial B({\bf 0},1)}  \Phi_{\omega}(x,y) \partial_{\nu} u_{k}(y)  d\sigma(y)  -  p.v. \int_{\partial B({\bf 0},1)}  \partial_{\nu}\Phi_{\omega}(x,y)  u_{k}(y)  d\sigma(y) =  0.
\end{equation*}
To solve \eqref{PT-Equa0918}, we start by splitting $u_{k}(\cdot): =  \hat{f_{k}}(\cdot)  +  \hat{f}_{k, H}(\cdot)$ in $B({\bf 0},1)$, where $\hat{f}_{k}(\cdot)$ satisfies
\begin{equation}\label{fk}
    \Delta \hat{f_{k}}(x)  +  \omega^{2}  \hat{f_{k}}(x) =  \, - \, f_{k}(x)\stackrel{(\ref{RBF})}{=} - \, \left( 1 + \left\vert x - x_{k} \right\vert \right), \quad x \in B({\bf 0},1),
\end{equation}
and $\hat{f}_{k, H}(\cdot)$ is the solution to the homogeneous Helmholtz equation
\begin{equation}\label{PT-Equa09180}
    \Delta \hat{f}_{k, H}(x)  +  \omega^{2}  \hat{f}_{k, H}(x)  =  0, \quad x \in B({\bf 0},1).
\end{equation}

Due to the fact that the right-hand side of \eqref{fk} is a real function, the possible imaginary part of $\hat{f}_{k}(\cdot)$ most satisfy the homogeneous Helmholtz equation. Therefore, we can add it to the function $\hat{f}_{k, H}(\cdot)$, the solution to \eqref{PT-Equa09180}. So we assume that $\hat{f}_{k}(\cdot)$ is a real function in the sequel, without lose of generality. Then, by plugging \eqref{PT-Equa0918} into \eqref{Equa0500}, and using \eqref{fk} and \eqref{PT-Equa09180}, we deduce
\begin{equation}\label{Eq0936II+}
    u_{k}(x)  =  \int_{B({\bf 0},1)}   \Phi_{\omega}(x,y) \left( \Delta \hat{f_{k}}(y)+\omega ^{2}\hat{f_{k}}(y) \right)  dy, \quad x \in B({\bf 0},1).
\end{equation}
Furthermore, in \eqref{Eq0936II+}, to handle the volume integral with weak singular kernel, we apply the Green's formula to change the volume integral into surface integral. Hence,
\begin{equation}\label{Eq0936II++}
    u_{k}(x)  =   \hat{f_{k}}(x) \, - \,  \int_{\partial B({\bf 0},1)} \left(  \Phi_{\omega}(x,y)  \frac{\partial \hat{f_{k}}}{\partial \nu(y)}(y) -  \hat{f_{k}}(y)  \frac{\partial \Phi_{\omega}}{\partial \nu(y)}(x,y) \right)  d\sigma(y), \quad x \in B({\bf 0},1).
\end{equation}
To compute numerically the surface integral in the right hand side of \eqref{Eq0936II++}, we need the analytic expression of $\hat{f_{k}}(\cdot)$. To keep things simple, we initially consider the case $x_{k} = 0$, and due to the Helmholtz equation being invariant when rotating, it's recommended to look for radial solution first, i.e., we take  $\hat{f}_{k}(x)  = g\left(\left\vert x \right\vert\right)$ for $x \in B({\bf 0},1)$. Then \eqref{fk} becomes
\begin{equation}\label{Eq1134}
    \left\vert x \right\vert  g^{\prime \prime}(\left\vert x \right\vert)  +  2  g^{\prime}(\left\vert x \right\vert)  +  \omega^{2}  \left\vert x \right\vert  g(\left\vert x \right\vert)  = \left\vert x \right\vert  +  \left\vert x \right\vert^{2}.
\end{equation}
By the function transform
$\varphi(\left\vert x \right\vert)  = \left\vert x \right\vert  g(\left\vert x \right\vert)$, the equation \eqref{Eq1134}  admits the  solution
\begin{equation*}\label{AB-Equa0751}
    \varphi(\left\vert x \right\vert)  =  \frac{\left\vert x \right\vert^{2}  +  \left\vert x \right\vert}{\omega^{2}}  -  \frac{2}{\omega^{4}} +  C_{1}  e^{-  i  \omega \left\vert x \right\vert}  +  C_{2}  e^{i  \omega  \left\vert x \right\vert} \quad \text{for} \quad C_1, C_2 \in \mathbb{C}.
\end{equation*}
This suggests, by the construction of $\varphi(\cdot)$, the solution to \eqref{fk} for general source point $x_k$ is
\begin{equation}\label{Equahatfk}
    \hat{f}_{k}(x)  =  \frac{\left\vert x - x_{k} \right\vert}{\omega^{2}} + \frac{1}{\omega^{2}} -  \frac{2}{\omega^{4}  \left\vert x - x_{k} \right\vert} + \frac{C_{1}  e^{-  i  \omega  \left\vert x - x_{k} \right\vert}}{\left\vert x - x_{k} \right\vert}  +  \frac{C_{2}  e^{i  \omega  \left\vert x - x_{k} \right\vert}}{\left\vert x - x_{k} \right\vert}, \quad x \in B({\bf 0},1),
\end{equation}
where $C_{1}$ and $C_{2}$ are such that
\begin{equation*}\label{C1C2}
    C_{1}   +  C_{2}   =  \frac{2}{\omega^{4}}.
\end{equation*}
To avoid complex solution to \eqref{Equahatfk}, we take $C_{1} =  C_{2}  = \omega^{-4}$ and then get
    \begin{equation}\label{hatf}
     \hat{f_{k}}(x)   =  \frac{\left( 1 +  \left\vert x - x_{k} \right\vert \right)}{\omega^{2}}  -  \frac{2 \left(1- \cos(\omega \left\vert x - x_{k} \right\vert)\right)}{\omega^{4}  \left\vert x - x_{k} \right\vert}.
    \end{equation}
Besides, by using the Taylor series related to the cosine function near 0 given by
\begin{equation*}
    \cos\left(\omega  \left\vert x - x_k \right\vert \right)  =  \sum_{n=0}^{\infty} \frac{(-1)^{n}}{(2n)!}  \left( \omega \left\vert x - x_k \right\vert \right)^{2n},
\end{equation*}
we deduce that $\hat{f_{k}}(x)\to \omega^{-2}$ as $\left\vert x - x_k \right\vert  \to  0$. Thus, $\hat{f_{k}}(\cdot)$ in $B({\bf 0},1)$ is of the removable singularity, which is a significant advantage for numerical computations. Now, by using  \eqref{fk} in  \eqref{Eq0936II}, we obtain for $x \in B({\bf 0},1)$ that
\begin{equation*}
 \sum_{k =  1}^{n} \alpha _{k} \left[ f_{k}(x)- \omega^{2}\left( \frac{1}{k_{0}(x)}-1 \right)
      \int_{B({\bf 0},1)}  \Phi_{\omega}(x,y) (\Delta \hat{f_{k}}(y)+\omega ^{2} \hat{f_{k}}(y))  dy \right] =  u^{\star,inc}(x, \theta, \omega)
        \end{equation*}
which, by an integration by parts, yields
        \begin{equation}\label{dis}
  \sum_{k =  1}^{n}\alpha _{k} \left[ f_{k}(x)-\omega^{2} \left( \frac{1}{k_{0}(x)}-1 \right)\left(\mathcal{J}_k(x,\omega)-\hat{f_{k}}(x)\right) \right]  =  u^{\star,inc}(x, \theta, \omega), \; x \in B({\bf 0},1),
        \end{equation}
 where
        \begin{equation*}\label{J}
            \mathcal{J}_k(x,\omega) := \int_{\partial B({\bf 0},1)} \left[  \Phi_{\omega}(x,y)\frac{\partial\hat{f_{k}}(y) }{\partial \nu(y)}  - \frac{\partial \Phi_{\omega}(x,y)}{\partial \nu(y)}\hat{f_{k}}(y) \right]  d\sigma(y), \quad x \in B({\bf 0},1).
        \end{equation*}
 Notice, for $\partial B({\bf 0},1)=\mathbb{S}^2$, the single and double surface potentials with known $\hat f_k(\cdot)$ given by  \eqref{hatf} involved in $\mathcal{J}_k(x,\omega)$ can be computed efficiently using the spherical harmonic functions expansions in $\mathbb{S}^2$, see  \cite{colton2019inverse} for more details. For any $x\in B({\bf 0},1)$, the integrand
 \begin{equation*}
      \mathcal{I}_k(y,x,\omega) :=  \Phi_{\omega}(x,y)\frac{\partial\hat{f_{k}}(y) }{\partial \nu(y)}  - \frac{\partial \Phi_{\omega}(x,y)}{\partial \nu(y)}\hat{f_{k}}(y),  \quad y \in \partial B({\bf 0},1)=\mathbb{S}^2,
 \end{equation*}
is an analytic function $\mathbb{S}^2 \to \mathbb{C}$, thus we have the Gauss trapezoidal product rule
\begin{equation}\label{J1}
    \mathcal{J}_k(x,\omega)\approx  \frac{\pi}{N}\sum_{j=1}^{N}\sum_{m=0}^{2N-1}\mu_j \mathcal{I}_k(y_{jm},x,\omega),
\end{equation}
where $\mu_j := \frac{2(1-t_j^2)}{[NP_{N-1}(t_j)]^2}$ for $ j= 1,\cdots,N$
are the weights of the Gauss-Legendre quadrature rule, with
$ -1 < t_{1} < t_{2} < \cdots < t_{N} < 1$ the zeros of the Legendre polynomial $P_{N}(\cdot)$, while $y_{jm}\in \mathbb{S}^2 $ are given in polar coordinates by
$y_{jm}:=(\sin \theta_j \cos \phi_m, \sin \theta_j \sin \phi_m, \cos \theta_j)$,
where $\theta_{j} = \arccos t_j$ for $j=1,\cdots,N$ and $\phi_m =\pi m/N$ for $m=0,\cdots,2N-1$.
Thus, by taking $x$ in \eqref{dis} as a random distributed collocation point $\left\{ x_{i} \right\}_{i=1}^{n} \subset B({\bf 0},1)$, we obtain
\begin{equation*}\label{xidis}
  \sum_{k =  1}^{n}\alpha _{k} \left[ f_{k}(x_{i})-\omega^{2}  \left( \frac{1}{k_{0}(x_i)}-1 \right)\left(\mathcal{J}_k(x_{i},\omega) - \hat{f_{k}}(x_{i})\right) \right] = u^{\star,inc}(x_i, \theta, \omega) \quad \text{for} \quad i=1,\cdots,n.
\end{equation*}
%for $i=1,\cdots,n$.
Then we derive the following algebraic system
\begin{equation}\label{alg}
    \mathbf{A} \cdot \mathbf{\alpha} =  \mathbf{b},
\end{equation}
where the elements of the interpolation matrix $\mathbf{A}=(A_{ik})\in \mathbb{C}^{n\times n }$, are given by
    \begin{equation*}\label{Aik}
    A_{ik} :=  f_{k}(x_{i})-\omega^{2}  \left( \frac{1}{k_{0}(x_i)}-1 \right)\left(\mathcal{J}_k(x_{i},\omega) - \hat{f_{k}}(x_{i})\right),  \quad \text{for} \,\, 1 \leq i, k \leq n,
\end{equation*}
the elements of the vector $\mathbf{b}=(b_1,\cdots,b_n)^T\in \mathbb{C}^{n\times 1}$, are given by $b_i=u^{\star,inc}(x_i, \theta, \omega)$, for $1 \leq i \leq n$,  and $\mathbf{\alpha}$ is the unknown vector give by  $\mathbf{\alpha}=(\alpha_1,\cdots,\alpha_n)^T\in \mathbb{C}^{n}$.

At first glance, one might say that the conditional positive definiteness of $f_{k}(x)$ renders the matrix $\mathbf{A}$ in \eqref{alg} potentially ill-conditioned for large $n$, which in turn compromises the invertibility of the algebraic system $\eqref{alg}$. Nevertheless, since $\mathbf{A}$ is composed of both the RBF terms $f_{k}(x_{i})$, with $1 \leq i, k \leq n$, and the frequency-dependent physical term $\mathcal{J}_k(x_i, \omega)$, with $1 \leq i, k \leq n$, the ill-posedness can be mitigated, thereby restoring stable invertibility, either by restricting to a moderate frequency $\omega$ (as in our case) or by incorporating a suitable regularization scheme.
%Since $f_{k}(x)$  is only conditionally positive definite, the matrix $\mathbf{A}$ in \eqref{alg} may be ill-conditioned for large $n$. However, the matrix $\mathbf{A}$ depends both on the RBF term $f_k(x_i)$ and on the frequency-dependent physical term $\mathcal{J}_k(x_i,\omega)$, we can improve the ill-posedness of $\mathbf{A}$ by considering moderate frequency $\omega$ or introduce  some regularization scheme directly.}

Finally, by solving   \eqref{alg}, we get the expansion coefficients $\alpha_{1},\cdots,\alpha_{n}$. Hence, the function $v^{\star}(\cdot, \theta, \omega)$ can be simulated in $B({\bf 0},1)$ using the expression \eqref{Vexpres}. Therefore, since we have assumed  that $k_{0}(\cdot)$ is {\it a-prior} known function, we deduce the reconstruction of $v(\cdot, \theta, \omega)|_{B({\bf 0},1)}$ by using \eqref{defvstar}. Consequently, from \eqref{Eq0823}, we deduce the reconstruction of $v(\cdot, \theta, \omega)$ in $\mathbb{R}^{3}$.

2. Computation of the total field $v^{\infty}(\cdot, \theta, \omega)$.

 Since $v(\cdot, \theta, \omega)  =  v^{s}(\cdot, \theta, \omega) +  u^{inc}(\cdot, \theta, \omega)$, the equation \eqref{Eq0823} says
    \begin{equation*}
        v^{s}(x, \theta, \omega)  =  \omega^{2} \int_{B({\bf 0},1)} \left( \frac{1}{k_{0}(y)}-1 \right)  \Phi_{\omega}(x,y) v(y, \theta, \omega) dy  \overset{(\ref{defvstar})}{=}
         \omega^{2}  \int_{B({\bf 0},1)} \Phi_{\omega}(x,y) v^{\star}(y, \theta, \omega)  dy, \; x \in \mathbb{R}^3.
    \end{equation*}
    Thanks to \cite[Theorem 2.6]{colton2019inverse}, it follows that $v^{\infty}(\cdot, \theta, \omega)$ has the following expression
    \begin{eqnarray}\label{V_infinity}
       v^{\infty}(\hat{x}, \theta, \omega)  = \frac{\omega^{2}}{4\pi}  \int_{B({\bf 0},1)}  e^{- i \omega \hat{x}  \cdot y}  v^{\star}(y, \theta, \omega) \, \overset{(\ref{Vexpres})}{=} \, \frac{\omega^{2}}{4\pi}  \sum_{k=1}^{n} \alpha_{k} \int_{B({\bf 0},1)} e^{-i \omega \hat{x} \cdot y} f_{k}(y)  dy, \quad \hat{x}\in \mathbb{S}^{2},
  \end{eqnarray}
 % which, by using $$, gives us
 % \begin{eqnarray}\label{V_infinity}
  %     v^{\infty}(\hat{x}, \theta, \omega) \, = \,
  %      \quad \hat{x}\in \mathbb{S}^{2},
 % \end{eqnarray}
where $f_{k}(\cdot)$ is given by \eqref{RBF}. That is, we can compute $v^{\infty}(\cdot, \theta, \omega)$ using $(\alpha_1,\cdots,\alpha_n)$ directly. In particular, by taking $\hat{x}  = -\theta$ in  \eqref{V_infinity}, we obtain
\begin{equation}\label{MS1036}
       v^{\infty}(-\theta, \theta, \omega)  =
       \frac{\omega^{2}}{4\pi}  \sum_{k=1}^{n} \alpha_{k} \int_{B({\bf 0},1)} e^{i \omega  \theta \cdot y}f_{k}(y)dy,
  \end{equation}
which can be computed numerically.

\subsection{Generated fields for the perturbed medium}\label{SubSecIII-III}

Similarly to the unperturbed medium case stated in the above subsection, we also generate the scattered wave related to the perturbed medium by the presence of the droplet $D_{z}$ in two parts. However, since the domain $\Omega \equiv B({\bf 0},1)$ involves the presence of a droplet $D_z$ with small size of order $\varepsilon^3$ and high contrasting bulk coefficient, see \eqref{defrhodefk} and \eqref{Ctek1}, the treatment for the volume integral should be careful.

In the presence of an injected droplet $D_{z}$, the total wave $u(\cdot,\cdot,\cdot)$ solving \eqref{Eq1} satisfies the following Lippmann-Schwinger equation
    \begin{equation}\label{LSEUR3}
        u\left( x, \theta, \omega \right)  -  \omega^{2}  \int_{B({\bf 0},1)}  \left( \frac{1}{k(y)}  -  1 \right)  \Phi_{\omega}(x,y)  u(y, \theta, \omega)  dy  =  u^{inc}(x, \theta, \omega), \quad x \in \mathbb{R}^{3},
    \end{equation}
    where the bulk modulus function $k(\cdot)$ is defined by \eqref{defrhodefk}. Notice, the solution to \eqref{LSEUR3} depends on $z$, which is denoted by $u_z(\cdot,\cdot,\cdot)$, since $k(\cdot)$ in $B({\bf 0},1)$ given by \eqref{defrhodefk} depends on $D_z$. In the same way as the unperturbed medium cases, we split the study into two steps.

1.  Estimation of the total field $u_z(\cdot,\cdot,\cdot)$.

By restricting \eqref{LSEUR3} into the unit ball $B({\bf 0},1)$, using \eqref{defrhodefk} and the fact that $\overline{k}_1=1$, we obtain
 \begin{eqnarray}\label{MLSA}
        u_z\left( x, \theta, \omega \right)  &-&  \omega^{2}  \int_{B({\bf 0},1) \setminus D_z}  \left( \frac{1}{k_{0}(y)}  -  1 \right)  \Phi_{\omega}(x,y)  u_z(y, \theta, \omega)  dy- \nonumber\\ & &  \omega^{2}  \left( \frac{1}{\varepsilon^{2}}  -  1 \right)  \int_{D_z}   \Phi_{\omega}(x,y)  u_z(y, \theta, \omega)  dy  =  u^{inc}(x, \theta, \omega), \quad x \in B({\bf 0},1).
    \end{eqnarray}
  To deal with the second term on the left hand side of \eqref{MLSA} which is an integral over a domain with a hole $D_z$, where the quadrature process is highly complicated  and may cause numerical instability, see \cite{John2009}, we rewrite \eqref{MLSA} as
    \begin{eqnarray}\label{comp_u}
     \nonumber
        u_z\left( x, \theta, \omega\right)  &-&  \omega^{2}  \int_{B({\bf 0},1)}   \Phi_{\omega}(x,y)  \left( \frac{1}{k_{0}(y)}  -  1 \right)  u_z(y, \theta, \omega)dy+ \\
        & & \omega^{2}  \int_{D_z}   \Phi_{\omega}(x,y) \left( \frac{1}{k_{0}(y)}  - \frac{1}{\varepsilon^{2}}  \right)  u_z(y, \theta, \omega) dy  =  u^{inc}(x, \theta, \omega), \quad x \in B({\bf 0},1).
    \end{eqnarray}

We solve the above equation using DRM directly. To do this, we start by approximating the function $u_z(\cdot,\theta,\omega)$, the solution to \eqref{comp_u},  by
\begin{equation}\label{liu011}
u_z(x,\theta,\omega)  =  \sum_{k=1}^n\beta_k f_k(x),\quad x\in B({\bf 0}, 1),
\end{equation}
where $f_k(\cdot)$ is the radial basis function defined by \eqref{RBF}. For notations' simplicity, we  omit the dependency of $\beta_{k} = \beta_{k}(z,\theta, \omega)\in \mathbb{C}$ on $(z,\theta,\omega)$.
By plugging \eqref{liu011} into \eqref{comp_u}, we have
    \begin{eqnarray*}
        \sum_{k=1}^n \beta_k f_k(x)  &-&  \omega^{2} \sum_{k=1}^n\beta_k  \int_{B({\bf 0},1)}   \Phi_{\omega}(x,y)  \left( \frac{1}{k_{0}(y)}  -  1 \right) f_k(y)  dy + \\
        & & \omega^{2}  \sum_{k=1}^n \beta_k  \int_{D_z}   \Phi_{\omega}(x,y) \left( \frac{1}{k_{0}(y)}  - \frac{1}{\varepsilon^{2}}  \right)  f_{k}(y) dy  =  u^{inc}(x, \theta, \omega), \quad x \in B({\bf 0},1),
    \end{eqnarray*}
which can be rewritten as
\begin{equation}\label{Equa0249}
    \sum_{k=1}^{n}  A_{k}(x)  \beta_{k}  =  u^{inc}(x, \theta, \omega), \quad x \in B({\bf 0},1),
\end{equation}
where
    \begin{eqnarray}\label{DefAk(x)}
    \nonumber
        A_{k}(x) &:=&  f_k(x)  -  \omega^{2} \int_{B({\bf 0},1)}   \Phi_{\omega}(x,y)  \left( \frac{1}{k_{0}(y)}  -  1 \right) f_k(y) dy \\ &+&   \omega^{2}  \int_{D_z}   \Phi_{\omega}(x,y) \left( \frac{1}{k_{0}(y)}  - \frac{1}{\varepsilon^{2}}  \right) f_{k}(y) dy.
    \end{eqnarray}
The evaluations of \eqref{Equa0249} at collocation points $\left\{ x_{1}, \cdots, x_{n} \right\}$ leads to the following algebraic system
\begin{equation}\label{A.beta=uinc}
    \begin{pmatrix}
        A_{1}(x_{1}) & \cdots & A_{n}(x_{1}) \\
        \vdots & \ddots & \vdots \\
        A_{1}(x_{n}) & \cdots & A_{n}(x_{n})
    \end{pmatrix} \cdot \begin{pmatrix}
        \beta_{1} \\
        \vdots \\
        \beta_{n}
    \end{pmatrix}  =  \begin{pmatrix}
        u^{inc}(x_{1}, \theta, \omega) \\
        \vdots \\
        u^{inc}(x_{n}, \theta, \omega)
    \end{pmatrix}.
\end{equation}

To generate the interpolation matrix appearing on the left hand side of \eqref{A.beta=uinc}, it is necessary to numerically compute the function $A_{k}(\cdot)$ given by \eqref{DefAk(x)} at the collocation points $\left\{ x_{1}, \cdots, x_{n} \right\}$. To do this, by using the fact that the bulk modulus function $k_{0}(\cdot)$, which is assumed {\it a-prior} known in simulation process, is a smooth function in $B({\bf 0},1)$, we have the following expansion
\begin{equation}\label{Equa327}
    \left( \frac{1}{k_{0}(y)}  - \frac{1}{\varepsilon^{2}}  \right)  f_{k}(y)  = \sum_{\ell = 1}^{n} c_{\ell,k}^{(\varepsilon)}  f_{\ell}(y), \quad y \in B({\bf 0},1),
\end{equation}
where $f_{\ell}(\cdot)$ is the radial basis function given by \eqref{RBF}. Besides, for known $k_{0}(\cdot)$ in $B({\bf 0},1)$, the coefficients $(c_{1,k}^{(\varepsilon)}, \cdots, c_{n,k}^{(\varepsilon)})$ are also known, which can be expressed by
$$c_{l,k}^{(\varepsilon)}=
\begin{cases}
d_{l,k}, &l\not=k, \\
d_{l,k}-\varepsilon^{-2}, &l=k,
\end{cases}
$$
where $\left( d_{1,k};\cdots;d_{n,k}\right)$ is solution of
\begin{equation*}\label{JM1012}
       \begin{pmatrix}
       f_{1}(x_{1}) & \cdots & f_{n}(x_{1}) \\
       \vdots & \ddots & \vdots \\
       f_{1}(x_{n}) & \cdots & f_{n}(x_{n})
    \end{pmatrix} \cdot \begin{pmatrix}
        d_{1,k}\\
        \vdots \\
        d_{n,k}
    \end{pmatrix}  =  \begin{pmatrix}
         \dfrac{1}{k_{0}(x_{1})}   f_{k}(x_{1}) \\
        \vdots \\
         \dfrac{1}{k_{0}(x_{n})}   f_{k}(x_{n})
    \end{pmatrix},
\end{equation*}
by taking the points $y=x_j\in B({\bf 0},1)$ in \eqref{Equa327} for $j=1,\cdots, n$.
It's worth noting that the component $c_{k,k}^{(\varepsilon)}$ is of a sharp jump $\varepsilon^{-2}$ for $0  <  \varepsilon \ll 1$, compared with $c_{l,k}^{(\varepsilon)}$ for $l\not=k$.

 Now, by plugging \eqref{Equa327} into \eqref{DefAk(x)}, we deduce for
 $x \in B({\bf 0},1)$ that
    \begin{equation}\label{EM-Equa1204}
        A_{k}(x)  =   f_k(x)  -  \omega^{2}   \sum_{\ell = 1}^{n} c_{\ell,k}^{(1)}  \int_{B({\bf 0},1)}   \Phi_{\omega}(x,y)  f_{\ell}(y)  dy
        +  \omega^{2} \sum_{\ell = 1}^{n} c_{\ell,k}^{(\varepsilon)}  \int_{D_z} \Phi_{\omega}(x,y)  f_{\ell}(y) dy,
    \end{equation}
where  $(c_{1,k}^{(1)}, \cdots, c_{n,k}^{(1)})$ are the expansion coefficients, such that
\begin{equation*}
    \left( \frac{1}{k_{0}(y)}  - 1  \right)  f_{k}(y)  = \sum_{\ell = 1}^{n} c_{\ell,k}^{(1)}  f_{\ell}(y), \quad y \in B({\bf 0},1).
\end{equation*}
By using \eqref{fk}, the equation \eqref{EM-Equa1204} can be written as
        \begin{eqnarray*}
        A_{k}(x)  &=&  f_k(x)   -  \omega^{2}  \sum_{\ell = 1}^{n} c_{\ell,k}^{(1)}  \int_{B({\bf 0},1)}   \Phi_{\omega}(x,y)   \left( \Delta \hat{f_{\ell}}(y)  +  \omega^{2} \hat{f_{\ell}}(y) \right)   dy+ \\
        & & \omega^{2}  \sum_{\ell = 1}^{n} c_{\ell,k}^{(\varepsilon)}  \int_{D_z}   \Phi_{\omega}(x,y)  \left( \Delta \hat{f_{\ell}}(y)  + \omega^{2}  \hat{f_{\ell}}(y) \right)   dy, \quad x \in B({\bf 0},1).
    \end{eqnarray*}

Thanks to Green's formula, the above equation becomes
\begin{eqnarray}\label{Equa0424}
        A_{k}(x)  &=&   f_k(x)  +  \omega^{2}   \sum_{\ell = 1}^{n} c_{\ell,k}^{(1)}  \hat{f_{\ell}}(x) -  \omega^{2}  \sum_{\ell = 1}^{n} c_{\ell,k}^{(\varepsilon)}  \hat{f_{\ell}}(x) \chi_{D_{z}}(x)-   \nonumber \\
        & & \omega^{2}  \sum_{\ell = 1}^{n} c_{\ell,k}^{(1)}   \int_{\partial B({\bf 0},1)} \left[\Phi_{\omega}(x,y)  \frac{\partial \hat{f_{\ell}}(y)}{\partial \nu(y)}  -  \frac{\partial \Phi_{\omega}(x,y)}{\partial \nu(y)}  \hat{f_{\ell}}(y)\right]  d\sigma(y) +\nonumber \\
        & & \omega^{2}  \sum_{\ell = 1}^{n} c_{\ell,k}^{(\varepsilon)}  \int_{\partial D_{z}} \left[  \Phi_{\omega}(x,y)  \frac{\partial \hat{f_{\ell}}(y)}{\partial \nu(y)}  - \frac{\partial \Phi_{\omega}(x,y)}{\partial \nu(y)} \hat{f_{\ell}}(y)\right]  d\sigma(y), \quad x \in B({\bf 0},1).
    \end{eqnarray}

For $x\in B({\bf 0},1)\setminus \partial D_z$, we can compute the single and double potentials defined on $\partial D_z$ in the above expression of $A_{k}(\cdot)$ using the formulas developed in \cite{Liu}. As for $x\in\partial D_z$, we can take some $\tilde x\in B({\bf 0},1)\setminus\overline{D}_z$ such that $\left\vert \tilde{x} -  x \right\vert  \to  0$. Then we have
\begin{eqnarray*}
    \int_{\partial D_{z}} \left[  \Phi_{\omega}(x,y)  \frac{\partial \hat{f_{\ell}}(y)}{\partial \nu(y)}  - \frac{\partial \Phi_{\omega}(x,y)}{\partial \nu(y)} \hat{f_{\ell}}(y)\right]  d\sigma(y)  &\approx &
\int_{\partial D_{z}} \left[  \Phi_{\omega}(\tilde x,y)  \frac{\partial \hat{f_{\ell}}(y)}{\partial \nu(y)}  - \frac{\partial \Phi_{\omega}(\tilde x,y)}{\partial \nu(y)} \hat{f_{\ell}}(y)\right]  d\sigma(y) \\&-& \frac{1}{2}\hat{f_{\ell}}(x), \quad x\in\partial D_z,
\end{eqnarray*}
due to the jump relation of the double layer potential, see \cite[Theorem 3.1]{colton2019inverse}. The integrals in the right hand side can also be computed in terms of the formulas developed in \cite{Liu}.
Finally, the matrix appearing on the left hand side of \eqref{A.beta=uinc} can be numerically computed for a set of collocation points $\{x_j:\; j=1,\cdots, n\}\subset B({\bf 0}, 1)$ with specified point $z\in B({\bf 0}, 1)$ determining the center of $D_z$. Solving \eqref{A.beta=uinc} enables us to determine the expansion coefficients $( \beta_{1}, \cdots, \beta_{n})$, and using \eqref{liu011} allows us to reconstruct the values of $u_{z}(\cdot, \theta, \omega)$, with $u_{z}(\cdot, \theta, \omega)$ being the solution to \eqref{LSEUR3}.

Furthermore, given the small size of $D_{z}$ and the limited number of collocation points, the expected number of points falling inside $D_{z}$ is effectively zero. In fact, for $\varepsilon = 0.01$ and $n=200$, the expected number of collocation points locating in $D_{z}$ is
\begin{equation*}
    E[D_{z}]  =  n  \frac{\left\vert D_{z} \right\vert}{\left\vert B({\bf 0}, 1) \right\vert}  =  n  \varepsilon^{3}  =  2 \times  10^{-4},
\end{equation*}
which is a value so close to zero that it is negligible in any practical context. Therefore, abstaining from the probabilistic interpretations associated with the above result, which fall outside the scope of our manuscript, we assume that all collocation points ${x_j:j=1,\cdots, n}$ are taken outside $D_z$, without loss of generality. As a result, the $\varepsilon^{-2}$-sharp jump for $x\in D_z$ originating from the third term in the expression of $A_k(x)$ (see the right-hand side of \eqref{Equa0424}) is circumvented.

In $(\ref{Equa0424})$, concerning the integral over the small surface $\partial D_z$, we give the following remark to provide further clarity.
\begin{remark}
In (\ref{Equa0424}), we need to compute integral
\[
I_\varepsilon(x):=\int_{\partial D_z}F(x,y)\,d\sigma(y),
\]
where
\[
F(x,y)=\Phi_\omega(x,y)\frac{\partial \hat f_\ell(y)}{\partial \nu(y)}
-\frac{\partial \Phi_\omega(x,y)}{\partial \nu(y)}\,\hat f_\ell(y).
\]
Since $\partial D_z=z+\varepsilon \mathbb{S}^2$,
under the transform $y=z+\varepsilon\hat y$ with $\hat y\in \mathbb{S}^2$, we have
\[
I_\varepsilon(x)=\varepsilon^2\int_{\mathbb{S}^2}G_\varepsilon(\hat y)\,d\sigma(\hat y),
\qquad
G_\varepsilon(\hat y):=F\bigl(x,z+\varepsilon\hat y\bigr).
\]
In our computation, the surface potentials on $\partial D_z$ are evaluated by the boundary element formulas in \cite[section 4.1]{Liu}.
More precisely, the discrete single-layer and double-layer potentials are computed by the coefficient matrices $Q_{ij}$ and $P_{ij}$ in \cite{Liu},
so the integral $I_\varepsilon(x)$ for any fixed $x$ is approximated by a weighted sum over the surface nodes.
We denote this discrete operator by $Q_N$ and write
\[
I_{\varepsilon,N}(x):=\varepsilon^2\,Q_N(G_\varepsilon)(x).
\]
Here $N$ denotes the surface discretization level such as the mesh size or the number of surface elements.
Define the quadrature error on $\mathbb{S}^2$ by
\[
E(N;g):=\left|\int_{\mathbb{S}^2} g(\hat y)\,d\sigma(\hat y)-Q_N(g)\right|.
\]
Then
\[
|I_\varepsilon(x)-I_{\varepsilon,N}(x)|=\varepsilon^2\,E(N;G_\varepsilon)(x).
\]
For $x\in B({\bf 0},1)\setminus D_z$, the integrand is smooth.
For $x\in\partial D_z$, we firstly use the jump relation stated after \eqref{Equa0424} to  handle the singular part,
and then treat the remaining smooth integrand. So we have
\[
E(N;G_\varepsilon)\le C\,E(N),
\]
where $E(N)\to0$ as the surface mesh is refined, and $C$ does not depend on $\varepsilon$.
Therefore,
\[
|I_\varepsilon(x)-I_{\varepsilon,N}(x)|\le C\,\varepsilon^2\,E(N).
\]
So for fixed $(x,N)$, the absolute quadrature error is $O(\varepsilon^2)$ as $\varepsilon\to 0$.
\end{remark}

2. Estimation of the far field $u_{z}^{\infty}\left(\cdot, \theta, \omega \right)$.

Similarly to the expression of $v^{\infty}(\cdot, \theta, \omega)$ given by \eqref{V_infinity}, the far-field $u_z^{\infty}(\cdot, \theta, \omega)$ associated to \\ $u_z^{s}(\cdot, \theta, \omega)  =  u_z(\cdot, \theta, \omega)-  u^{inc}(\cdot, \theta, \omega)$, with $u_z(\cdot, \theta, \omega)$  solving \eqref{LSEUR3}, is given by
    \begin{eqnarray*}
       u_z^{\infty}(\hat{x}, \theta, \omega)  &=& \omega^{2}\int_{B({\bf 0},1)}  \left( \frac{1}{k(y)} -1 \right)\frac{e^{- i \omega\hat{x} \cdot y}}{4 \pi}   u_z(y, \theta, \omega) dy \\
    &\stackrel{ (\ref{liu011})}{=}&
    \omega^{2} \sum_{k=1}^n\beta_k\int_{B({\bf 0},1)}  \left( \frac{1}{k(y)}  -  1 \right)  \frac{e^{- i  \omega \hat{x}  \cdot y}}{4  \pi}  f_k(y)  dy,
    \end{eqnarray*}
    where the expansion coefficients $(\beta_1,\cdots,\beta_n)$ are solved from \eqref{A.beta=uinc}. Besides, by using \eqref{defrhodefk}, \eqref{Ctek1}, and the face that $\overline{k}_1=1$, the above equation for $\hat{x} = -\theta$ becomes
    \begin{equation}\label{liu019}
    u_z^{\infty}(-\theta, \theta, \omega)   =   \frac{\omega^{2}}{4  \pi} \sum_{k=1}^n\beta_k \left[ \int_{B({\bf 0},1)} e^{i  \omega  \theta  \cdot y}  \left( \frac{1}{k_{0}(y)}  -  1 \right)    f_k(y)  dy-  \int_{D_{z}} e^{i  \omega  \theta  \cdot y}   \left( \frac{1}{k_{0}(y)}  - \frac{1}{\varepsilon^{2}}  \right)  f_k(y)  dy\right].
\end{equation}

The expression \eqref{MS1036} (respectively, \eqref{liu019}) construct an efficient scheme for direct scattering problem for an inhomogeneous medium without (respectively, with) droplet, which also provide us the simulated data for inverse scattering problem,
see \eqref{BSDI}. We will apply these data to check the efficiency of our reconstruction scheme stated in \textbf{Algorithm \ref{AlgoScheme}}.

    To conclude this section, we include the following clarification regarding the efficiency of our approximate scheme.
\begin{remark}
    Our developed computational algorithm for solving the forward scattering problem is quiet efficient, as compared with the straightforward computations from the Lippmann-Schwinger \eqref{liu04} for unperturbed medium, and \eqref{comp_u} for perturbed medium. Once we expand  the scattered wave in terms of the specified base functions with the expansion coefficients $\{\alpha_k: k=1,\cdots, n\}$, see  \eqref{Vexpres} for example,  the computational complexity of two schemes comes from the computational costs for all the elements of the coefficient matrix, if the corresponding linear systems are assumed to be solved by the same numerical scheme such as Gauss-Seidal iteration scheme. For $N$-sampling points $\{x_j:j=1,\cdots,N\}\subset B({\bf 0},1)$, the straightforward computations for the matrix elements from the Lippmann-Schwinger equation needs to compute $N^2$-volume integrals with weak singularity. However, by choosing special bases functions and applying the DRM scheme for solving \eqref{liu04}, which transfers the volume integrals into surface integrals, our scheme constitutes the matrix elements by $N^2$-surface integrals $\mathcal{J}_k(x_i,\omega)$ in $\partial B({\bf 0},1)$ with smooth integrands for $x\in B({\bf 0},1)$ in terms of \eqref{J1}, which are relatively cheap.
\end{remark}

%%%%%%%%%%%%%%%%%%%%%%%%%%%%%%%%%%%%%%%%%%%%%%%%%%%%%
%%%%%%%%%%%%%%%%%%%%%%%%%%%%%%%%%%%%%%%%%%%%%%%%%%%%%
\section{Numerical implementations for direct scattering}\label{Sec-Num-Sim}

This section is used to check the numerical implementations of our integral solver and the far-field computation. Here we use a manufactured-solution test. We choose a smooth complex-valued function, and then build the right-hand side so that the smooth complex-valued function satisfies the integral equation. So this test checks the code and the discretization. It is not a physical scattering test. Such a realization process also provides us the simulation data for the far fields patterns $u_z^{\infty}(-\theta,\theta,\omega)$ and $v^{\infty}(-\theta,\theta,\omega)$. Then by moving $D_z$ inside $\Omega$, we collect the data $\{u_z^\infty(-\theta,\theta,\omega):\ D_{z} \subset \subset \Omega\}$ and  $v^\infty(-\theta,\theta,\omega)$ from which we can recover $k_0(z)$ in $\Omega$, as explained in \textbf{Algorithm \ref{AlgoScheme}}. In the recovery process, we take $\varepsilon  =  0.01$, $h  =  0.95$, and we only need to use one fixed incident direction $\theta  = \dfrac{1}{\sqrt{6}}  \left(1, 2, 1\right) \in \mathbb{S}^2$.

In this section, we only test our algorithm for moderate frequency $\omega$ given on Table \ref{Table1}, with $n_0=1$.  It is well-known that the solutions to the Helmholtz equations with large frequencies oscillate seriously. The implementations of our scheme for higher frequencies $\omega$ are typically struggle, which should be studied in future research.

\subsection{The unperturbed medium case}\label{SubSectionIII-II}

To check our algorithm,  we recall from \eqref{liu05} that
 \begin{equation}\label{liu05+}
         v^{\star}(x, \theta, \omega)  -  \omega^{2}
        \left( \frac{1}{k_{0}(x)}  -  1 \right) \int_{B({\bf 0},1)}   \Phi_{\omega}(x,y)   v^{\star}(y, \theta, \omega)  dy  =  u^{\star,inc}\left(x, \theta, \omega\right), \quad x \in B({\bf 0},1).
    \end{equation}
The exact solution $v^{\star}(\cdot, \theta, \omega)$ for \eqref{liu05+} for general $k_0(\cdot)$, although exists, is very hard to be expressed analytically. Then, our goal is to obtain its approximation using our algorithm, which we refer to as $v^{\star}_{num}(\cdot, \theta, \omega)$. To verify the validity of our scheme, we replace the right hand side of \eqref{liu05+} by a known source function $\Theta(\cdot, \theta, \omega)$ generated from some assumed a priori known function $v^{\star}(\cdot,\theta, \omega)$ with analytical expression, that is,
\begin{equation}\label{DefTheta}
        \Theta(x, \theta, \omega) := v^{\star}(x, \theta, \omega)  -  \omega^{2}
        \left( \frac{1}{k_{0}(x)}  -  1 \right) \int_{B({\bf 0},1)}   \Phi_{\omega}(x,y) v^{\star}(y, \theta, \omega)  dy.
\end{equation}
Hence, in \eqref{liu05+}, by replacing the right hand side with $\Theta(\cdot, \theta, \omega)$ and using our algorithm, we obtain the numerical solution $v^{\star}_{num}(\cdot, \theta, \omega)$ associated to the known function  $v^{\star}(\cdot, \theta, \omega)$. The objective is to justify the small difference between $v^{\star}(\cdot, \theta, \omega)$, which is known, and $v^{\star}_{num}(\cdot, \theta, \omega)$ which will be computed using our numerical scheme, see \textbf{Algorithm \ref{AlgoScheme}}. To this end, we consider the case that $v^{\star}(\cdot,\theta, \omega)$ and $k_{0}(\cdot)$ are given by the following expressions
\begin{equation*}\label{DefvstarDefk0}
    v^{\star}(x, \theta, \omega) =  \left| x -x_1^*\right|^2  +  i  \left| x -x_2^*\right|^2, \quad \text{and} \quad k_0(x) = \frac{2}{1+\left| x \right|^2} \quad \text{for} \quad x\in B({\bf 0},1),
\end{equation*}
respectively, where $x_1^*=(1,-1.5,1.5)$ and $x_2^*=(1.5,0,1.5)$. Then \eqref{DefTheta} becomes
\begin{eqnarray}\label{FFJM}
   \Theta(x,\theta, \omega)  =  \frac{ \left| x \right|^2 + 1 }{2}  ( \left| x -x_1^*\right|^2+i\left| x -x_2^*\right|^2 ) -  \frac{( \left| x \right|^2 -  1 )  ( 3  +  i  3)}{\omega^{2}} +
        \frac{\left| x \right|^2 -  1 }{2}  \left[  \mathcal{F}_{1}(x) - \mathcal{F}_{2}(x) \right]
\end{eqnarray}
from the Green's formula, where the boundary integrals $\mathcal{F}_{1}(\cdot)$ and $\mathcal{F}_{2}(\cdot)$ are given by
\begin{eqnarray*}
   \mathcal{F}_{1}(x)  &:=&  \int_{\partial B({\bf 0},1)}    \frac{\partial \Phi_{\omega}(x,y)}{\partial \nu(y)}   \left( \left| y -x_1^*\right|^2+i\left| y -x_2^*\right|^2 -  \frac{(6+6i)}{\omega^{2}} \right) d\sigma(y), \quad x \in B({\bf 0},1), \\
   \mathcal{F}_{2}(x) &:=&  \int_{\partial B({\bf 0},1)}    \Phi_{\omega}(x,y)  \frac{\partial  \left( \left| y -x_1^*\right|^2+i\left| y -x_2^*\right|^2 \right)}{\partial \nu(y)} d\sigma(y), \qquad \qquad \quad   \quad x \in B({\bf 0},1),
\end{eqnarray*}
which can be calculated numerically with high accuracy like $\mathcal{J}_k(x,\omega)$, see \eqref{J1}. Hence, the source term $\Theta(\cdot, \theta, \omega)$ given by \eqref{FFJM} can be numerically computed. Consequently, by using our proposed algorithm scheme, the numerical solution $v_{num}^{\star}(\cdot,\cdot,\cdot)$ for $x \in B({\bf 0},1)$ associated to such a known function $v^{\star}(\cdot,\cdot,\cdot)$ can be derived by solving the following equation
 \begin{equation}\label{v_num}
         v_{num}^*(x, \theta, \omega)  -  \omega^{2}
        \left( \frac{1}{k_{0}(x)}  -  1 \right) \int_{B({\bf 0},1)}   \Phi_{\omega}(x,y)   v_{num}^*(y, \theta, \omega)  dy = \Theta(x,\theta,\omega),
    \end{equation}
where $\Theta(x,\theta,\omega)$ is given by \eqref{FFJM}.

\begin{figure}[H]
    \centering
    \subfigure[Comparison of the real part for a polar radius of $r=0.3$.]{
    \centering
    \includegraphics[width=1\linewidth,height=0.22\textwidth]{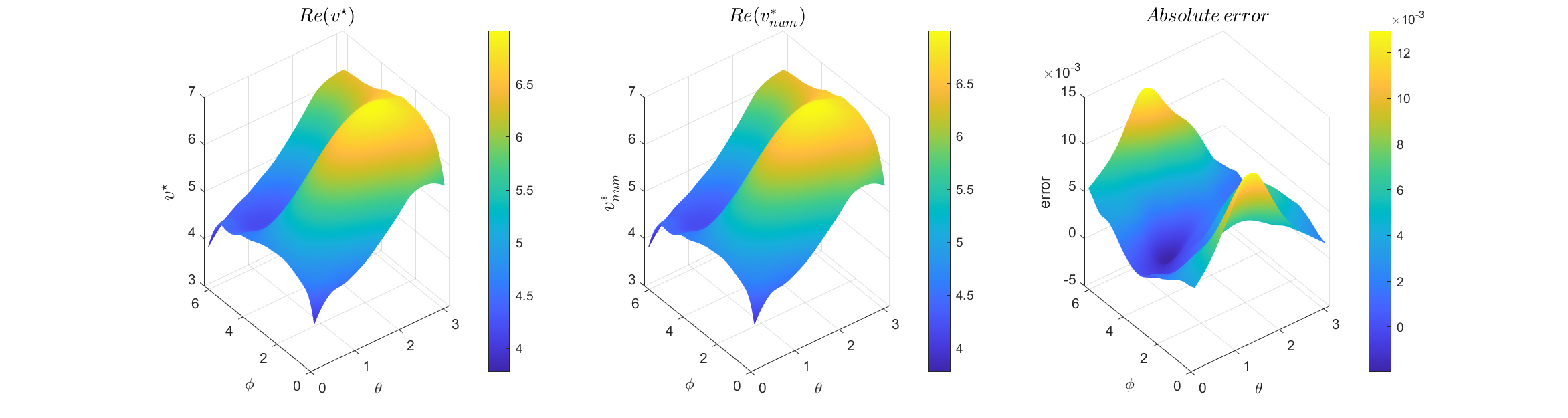}
    }
    \subfigure[Comparison of the imaginary part for a polar radius of $r=0.3$.]{
    \centering
    \includegraphics[width=1\linewidth,height=0.22\textwidth]{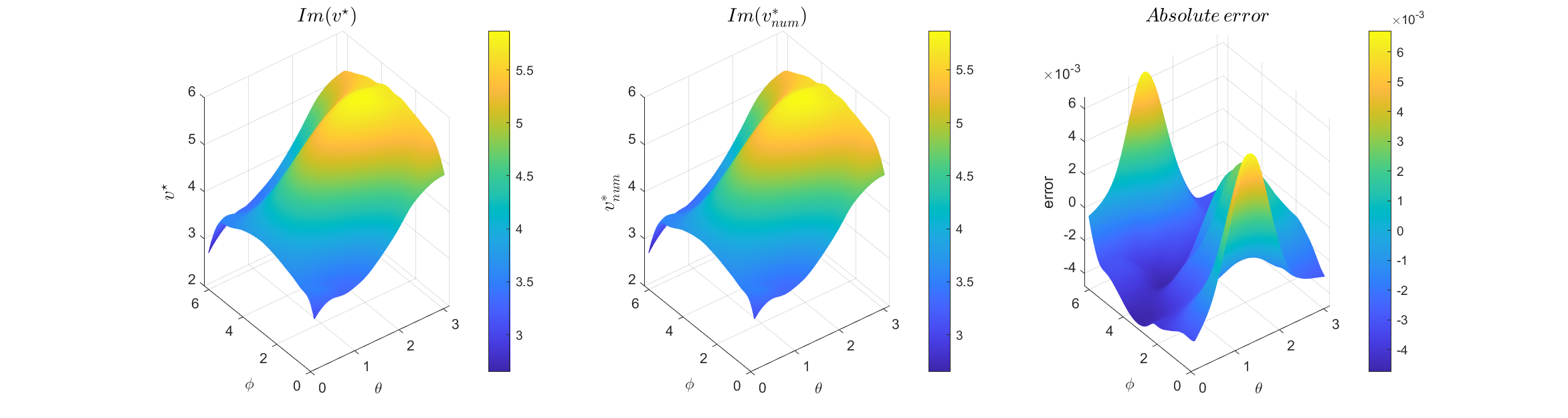}
    }
    \caption{Comparison of $\phi-\tilde{\theta}$ distributions at polar radius $r=0.3$.}
    \label{Comp03}
\end{figure}

\begin{figure}[h!]
    \centering
    \subfigure[Comparison of the real part for a polar radius of $r=0.6$.]{
    \centering
    \includegraphics[width=1\linewidth,height=0.22\textwidth]{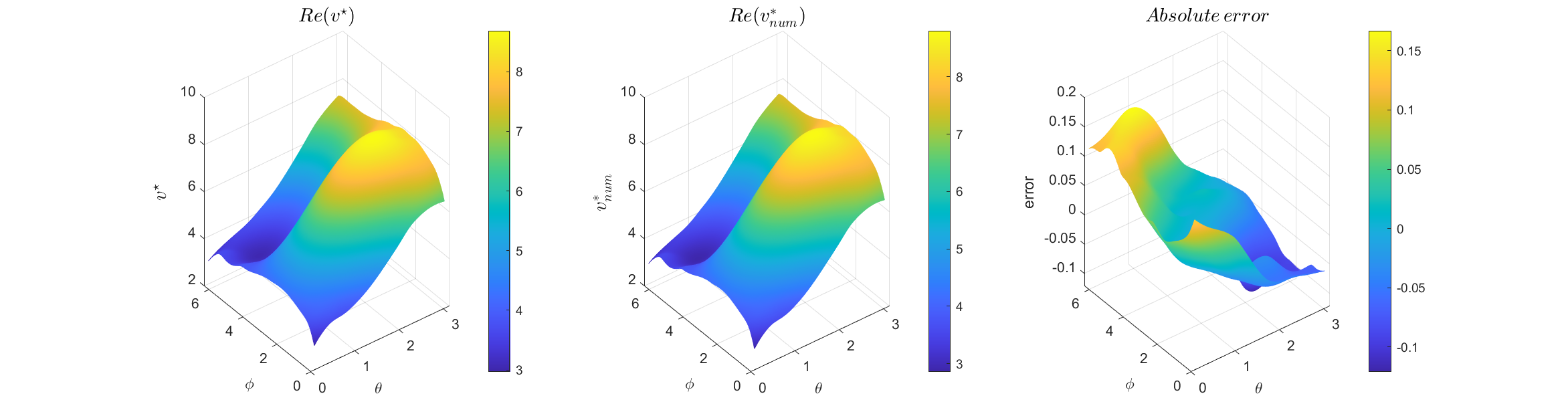}
    }
    \subfigure[Comparison of the imaginary part for a polar radius of $r=0.6$.]{
    \centering
    \includegraphics[width=1\linewidth,height=0.22\textwidth]{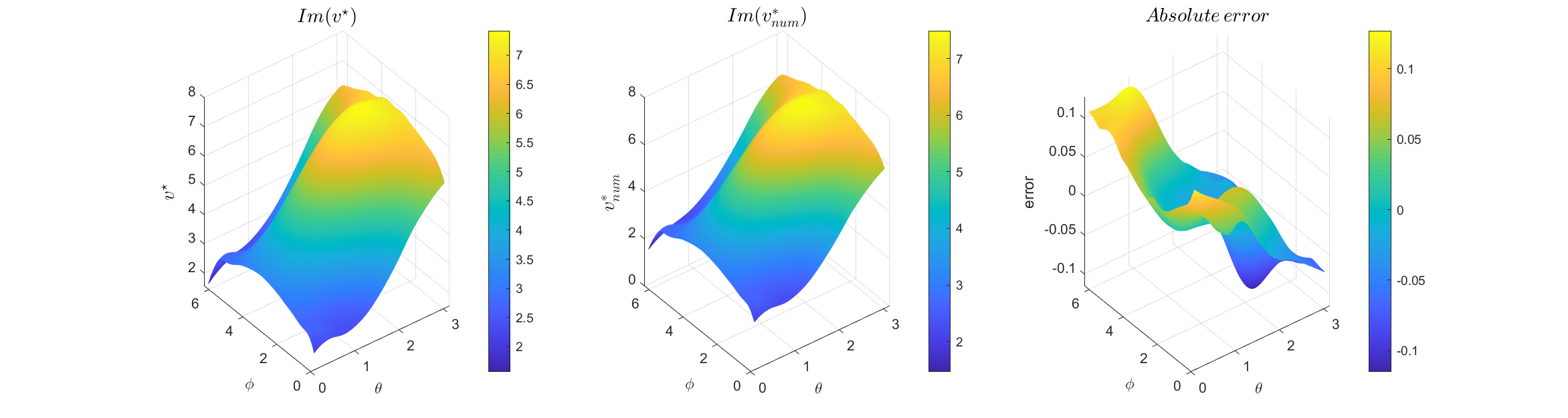}
    }
    \caption{Comparison of $\phi-\tilde{\theta}$ distributions at polar radius $r=0.6$.}
    \label{Comp06}
\end{figure}

\begin{figure}[h!]
    \centering
    \subfigure[Comparison of real part of the far-field pattern.]{
    \centering
    \includegraphics[width=1\linewidth,height=0.22\textwidth]{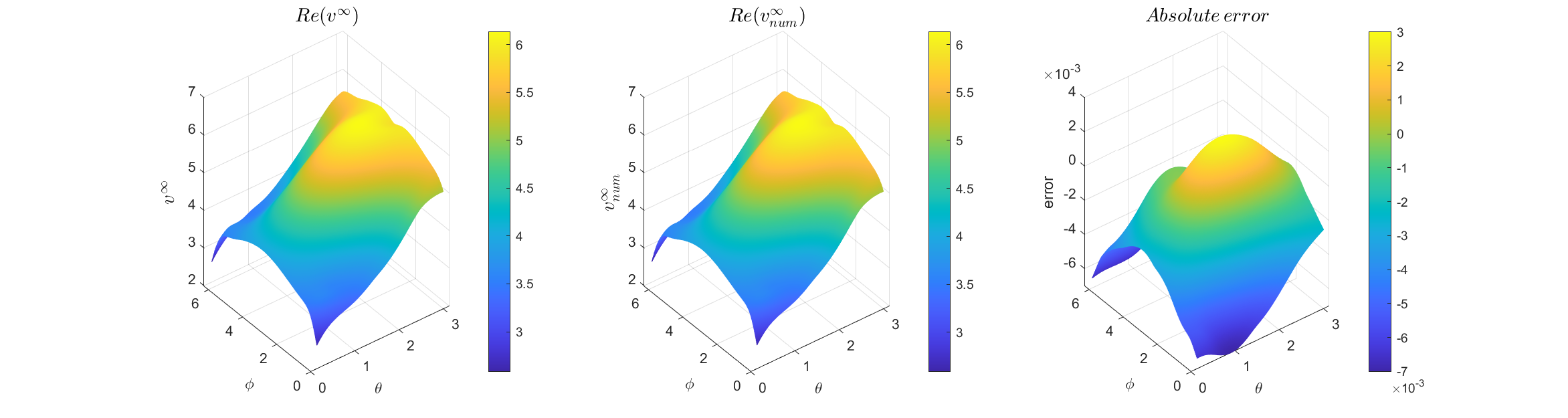}
    }
    \subfigure[Comparison of imaginary part of the far-field pattern.]{
    \centering
    \includegraphics[width=1\linewidth,height=0.22\textwidth]{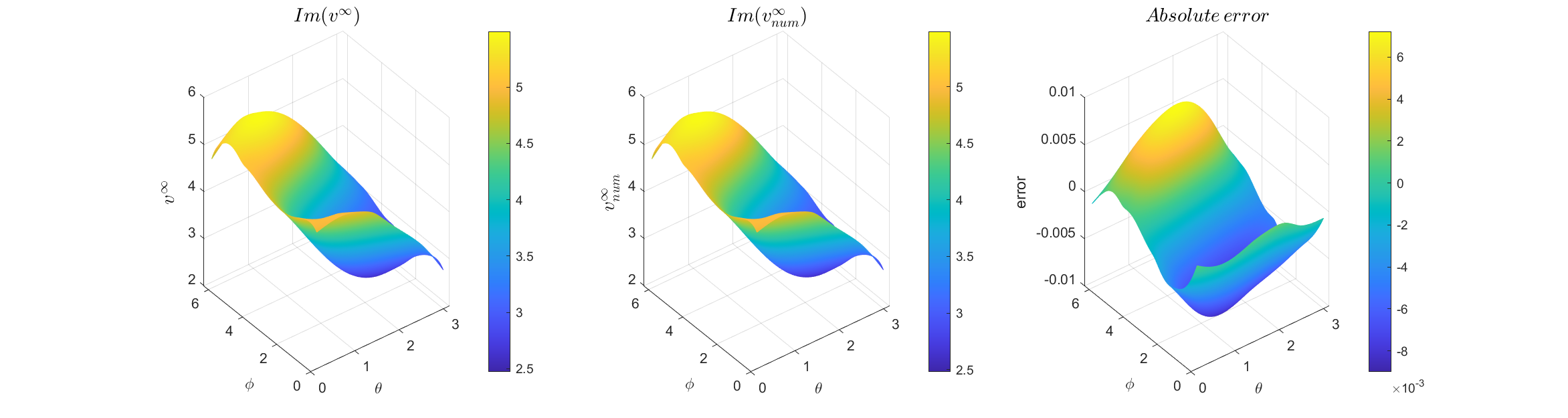}
    }
    \caption{Comparison of far-field for $\hat x\in \mathbb{S}^2$.}
    \label{Compvinf}
\end{figure}

The accuracy of our proposed algorithm scheme can be tested by comparing the disparity between the assumed known function $v^{\star}(\cdot,\theta,\omega)$ and its associated numerical solution $v_{num}^{\star}(\cdot,\theta,\omega)$, solution of \eqref{v_num}. In order to do this, for both $v^{\star}(x,\theta,\omega)$ and $v_{num}^{\star}(x,\theta,\omega)$, we represent the argument $x \in B({\bf 0},1) $ by spherical coordinates, i.e., $x =  r  \left( \sin\tilde{\theta} \cos\phi, \sin\tilde{\theta}\sin\phi, \cos\tilde{\theta} \right)$, with  $r \in (0,1), \tilde{\theta}\in (0,2\pi)$, and  $\phi \in (0,\pi)$, and we use spherical surfaces with polar radius of $r  =  0.3$, and $r  =  0.6$ in polar coordinates to compare between $v^{\star}(r,\tilde{\theta}, \phi,\theta,\omega)$ and $v_{num}^{\star}(r,\tilde{\theta}, \phi,\theta,\omega)$.  In Figure \ref{Comp03} and Figure \ref{Comp06},  we give the comparisons of real part and imaginary part between $v^{\star}(\cdot,\theta,\omega_{1})$ and $v^{\star}_{num}(\cdot,\theta,\omega_{1})$ for radius $r = 0.3,0.6$, respectively, while the comparison of real part and imaginary part of  $v^{\infty}(\cdot,\theta,\omega_{1})$ and $v^{\infty}_{num}(\cdot,\theta,\omega_{1})$ in terms of \eqref{V_infinity} for known expansion coefficients are shown in Figure \ref{Compvinf}. These three figures show the validity of our proposed scheme.

We use two different error measures to describe the numerical performances.
First, we use the pointwise absolute error on the test surfaces given by
\[
E_{\mathrm{pw}}^{\mathrm{re}}(x):=
\left|\Re\!\left(v^\star(x,\theta,\omega_1)-v^\star_{\mathrm{num}}(x,\theta,\omega_1)\right)\right|,
\quad
E_{\mathrm{pw}}^{\mathrm{im}}(x):=
\left|\Im\!\left(v^\star(x,\theta,\omega_1)-v^\star_{\mathrm{num}}(x,\theta,\omega_1)\right)\right|
\]
to measure the numerical performance of our forward solver for real part and imaginary one.
Figure \ref{Comp03} and Figure \ref{Comp06} show that both $E_{\mathrm{pw}}^{\mathrm{re}}$ and $E_{\mathrm{pw}}^{\mathrm{im}}$ are mostly of order $10^{-2}$ on the test surfaces.
  Second, we use the absolute $L^2$ error to indicate the numerical performance in the whole domain $B({\bf 0},1)$ given by
\[
E_2:=\left\|v^\star(\cdot,\theta,\omega_1)-v^\star_{\mathrm{num}}(\cdot,\theta,\omega_1)\right\|_{L^2(B({\bf 0},1))}
=\left(\int_{B({\bf 0},1)}\left|v^\star(x,\theta,\omega_1)-v^\star_{\mathrm{num}}(x,\theta,\omega_1)\right|^2dx\right)^{1/2}.
\]
This error quantifies numerical performances of real part and imaginary one together in the ball $B({\bf 0},1)$, not just in the surface $B({\bf 0},r_0)$ for $r_0\in (0,1)$. We report this absolute $L^{2}(B({\bf 0},1))$-norm in terms of \eqref{Vexpres} for computed expansion coefficients that
\begin{equation}\label{WZEq1}
E_2 \simeq 0.036.
\end{equation}

Obviously, our algorithm scheme exhibits excellent numerical performances of scattered wave without droplet  from the numerical experiments. Consequently, our proposed scheme can be utilized to obtain a numerical solution  $v_{num}^{\star}(\cdot, \theta, \phi)$ to \eqref{liu05+}, which approximates the unknown function $v^{\star}(\cdot, \theta, \phi)$. Hence, by solving \eqref{liu05+} with our proposed algorithm, we can determine the coefficients $(\alpha_{1}, \cdots, \alpha_{n})$. These coefficients also enable us to obtain an approximation to the numerical far-field $v^{\infty}_{num}(\hat{x}, \theta, \phi)$ related to $v_{num}(x, \theta, \phi)$ from (\ref{V_infinity}).
%, not only the good approximation to $v^{\star}(\cdot, \theta, \phi)$ using (\ref{Vexpres}).
The numerical approximations of the far-field pattern for $\hat x$ expressed by polar coordinates are shown in Figure \ref{Compvinf}.
It is apparent that the performance of $v^{\infty}_{num}(\hat{x}, \theta, \phi)$ is better than that related to the numerical solution $v_{num}(\hat{x}, \theta, \phi)$. This is quite natural, since it is well-known that $v^{\infty}_{num}(\hat{x}, \theta, \phi)$ is generated by integrating $v_{num}(\hat{x}, \theta, \phi)$, see for instance \eqref{V_infinity}, which acts as a regularizing operator in some sense by the average effect of the integration process.

%%%%%%%%%%%%%%%%%%%%%%%%%%%%%%%%%%%%%%%%%%%%%%%%%%%%%
%%%%%%%%%%%%%%%%%%%%%%%%%%%%%%%%%%%%%%%%%%%%%%%%%%%%%
\subsection{The perturbed medium case}

Now we test the simulation process for $u_z(\cdot,\theta, \omega)$ as well as $u_z^\infty(\cdot,\theta, \omega)$ numerically. To do this, we recall from \eqref{comp_u} that we have
\begin{eqnarray}\label{Equa0716}
     \nonumber
        u_z\left( x, \theta, \omega\right)  &-&  \omega^{2}  \int_{B({\bf 0},1)}   \Phi_{\omega}(x,y)  \left( \frac{1}{k_{0}(y)}  -  1 \right)  u_z(y, \theta, \omega)dy+ \\
        & & \omega^{2}  \int_{D_z}   \Phi_{\omega}(x,y) \left( \frac{1}{k_{0}(y)}  - \frac{1}{\varepsilon^{2}}  \right)  u_z(y, \theta, \omega) dy  =  u^{inc}(x, \theta, \omega), \quad x \in B({\bf 0},1).
    \end{eqnarray}
As done in \textbf{Subsection \ref{SubSectionIII-II}}, we firstly assume that the two functions $k_{0}(\cdot)$ and $u_z( \cdot, \theta, \omega)$ are known. Besides, we denote the source term $\Theta( \cdot, \theta, \omega)$ with known $u_z( \cdot, \theta, \omega)$ given by
\begin{eqnarray}\label{Equa0728}
     \nonumber
        \Theta\left( x, \theta, \omega\right)  &:=&  u_{z}\left( x, \theta, \omega\right)  -  \omega^{2}  \int_{B({\bf 0},1)}   \Phi_{\omega}(x,y)  \left( \frac{1}{k_{0}(y)}  -  1 \right)  u_{z}(y, \theta, \omega)dy+ \\
        & & \omega^{2}  \int_{D_z}   \Phi_{\omega}(x,y) \left( \frac{1}{k_{0}(y)}  - \frac{1}{\varepsilon^{2}}  \right)  u_{z}(y, \theta, \omega) dy, \quad x \in B({\bf 0},1).
    \end{eqnarray}
Then,  by replacing the right hand side in \eqref{Equa0716} with $\Theta( \cdot, \theta, \omega)$  and using our algorithm, we obtain the numerical solution $u_{z,num}( \cdot, \theta, \omega)$ associated to the known function $u_{z}( \cdot, \theta, \omega)$. To verify the efficiency of our proposed scheme for solving \eqref{Equa0716} numerically, we check the difference between known $u_{z}( \cdot, \theta, \omega)$ and $u_{z,num}( \cdot, \theta, \omega)$, solution of
\begin{eqnarray*}
     \nonumber
        u_{z,num}\left( x, \theta, \omega\right)  &-&  \omega^{2}  \int_{B({\bf 0},1)}   \Phi_{\omega}(x,y)  \left( \frac{1}{k_{0}(y)}  -  1 \right)  u_{z,num}(y, \theta, \omega)dy+ \\
        & & \omega^{2}  \int_{D_z}   \Phi_{\omega}(x,y) \left( \frac{1}{k_{0}(y)}  - \frac{1}{\varepsilon^{2}}  \right)  u_{z,num}(y, \theta, \omega) dy  =  \Theta\left( x, \theta, \omega\right), \; x \in B({\bf 0},1),
    \end{eqnarray*}
    where $\Theta\left( x, \theta, \omega\right)$ is given by \eqref{Equa0728}. We consider the configuration where $z = (0.2,0.3,0.6)$ and
\begin{equation*}\label{DefvstarDefk0per}
    u_{z}(x, \theta, \omega)  =  \left| x \right|^2+\textbf{b}\cdot x+3,\quad \text{and} \quad k_0(x) =  \frac{\left| x \right|^2  +  i \textbf{b}\cdot x  +3}{4\left| x \right|^2+i \, \textbf{b} \cdot x}, \quad \quad x\in B({\bf 0},1),
\end{equation*}
with $\textbf{b} = (1,2,3)$.
% Hence,  by replacing the right hand side of (\ref{Equa0716}) with $\Theta( \cdot, \theta, \omega)$ determined from (\ref{Equa0728}) and (\ref{DefvstarDefk0per}), we obtain the numerical approximation related to the solution $\textcolor{red}{u_{z}}( \cdot, \theta, \omega)$, i.e., $\textcolor{red}{u_{z,num}}( \cdot, \theta, \omega)$ is solved numerically from
%\begin{eqnarray*}
%        \textcolor{red}{u_{z,num}}\left( x, \theta, \omega\right)  &-&  \omega^{2}  \int_{B(0,1)}   \Phi_{\omega}(x,y)  \left( \frac{1}{k_{0}(y)}  -  1 \right)  \textcolor{red}{u_{z,num}}(y, \theta, \omega)dy \\
%        &+& \omega^{2}  \int_{D_z}   \Phi_{\omega}(x,y) \left( \frac{1}{k_{0}(y)}  - \frac{1}{\varepsilon^{2}}  \right)  \textcolor{red}{u_{z,num}}(y, \theta, \omega) dy  =  \Theta(x, \theta, \omega), \quad x \in B(0,1).
%    \end{eqnarray*}
    The accuracy of our proposed algorithm scheme can be tested by comparing the disparity between the assumed known function $u_{z}( \cdot, \theta, \omega)$ and its associated numerical solution $u_{z,num}(\cdot, \theta, \omega)$. In order to do this, for both $u_{z}( x, \theta, \omega)$ and $u_{z,num}(x, \theta, \omega)$, we represent the argument $x$  in terms of the spherical coordinate
    $x = r_0(1,0,0)+r
(\sin\tilde{\theta} \cos\phi; \sin\tilde{\theta}\sin\phi; \cos\tilde{\theta})$,
where the two parameters $r_{0}, r$ should be taken small such that $x\in B({\bf 0},1)$.  In our numerics, we fix $r_0=0.1$, and take $r = 0.2$ and $r =  0.4$ respectively, to compare
$u_{z}(r, \tilde{\theta}, \phi, \theta, \omega)$ with $u_{z,num}( r, \tilde{\theta}, \phi, \theta, \omega)$. The following figures present the comparisons in two spherical surfaces with radius $r=0.2$ and $r=0.4$.

\begin{figure}[h!]
    \centering
    \subfigure[Comparison of the real part for  polar radius  $r=0.2$.]{
    \centering
    \includegraphics[width=1\linewidth,height=0.22\textwidth]{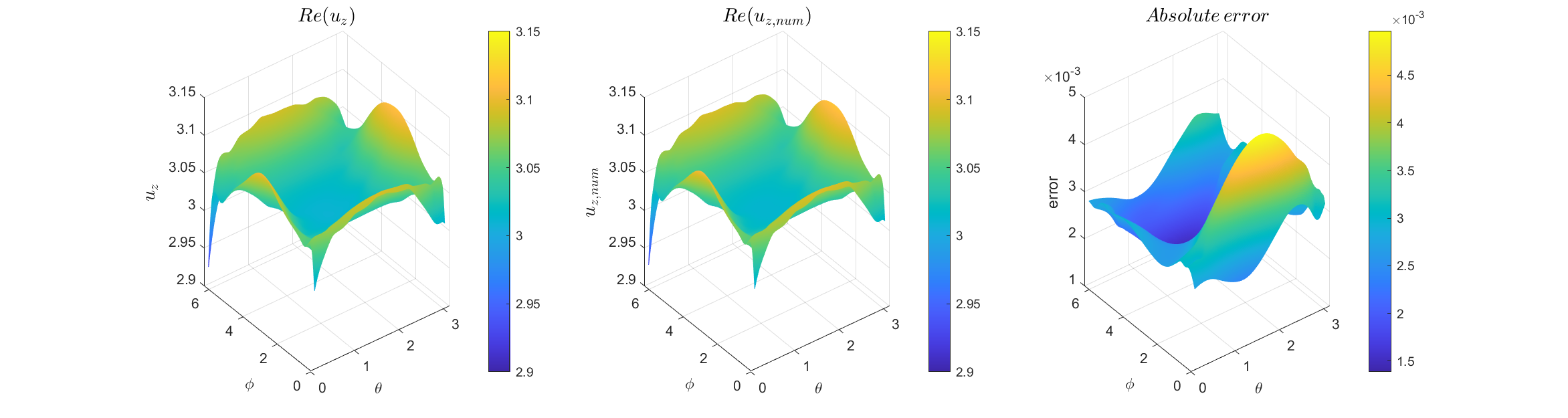}
    }
    \subfigure[Comparison of the imaginary part for  polar radius  $r=0.2$.]{
    \centering
    \includegraphics[width=1\linewidth,height=0.22\textwidth]{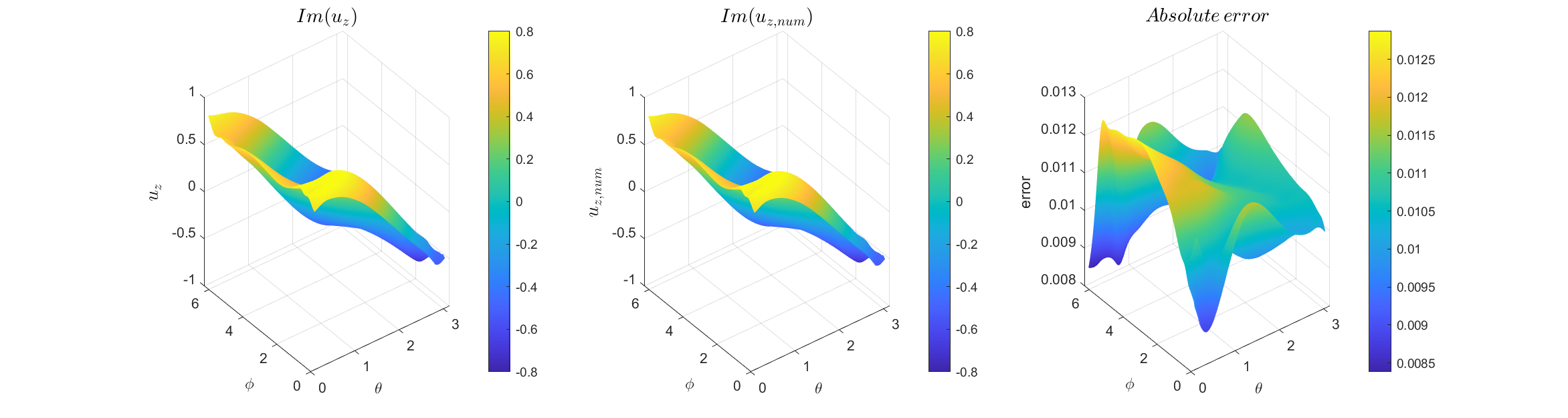}
    }
    \caption{Comparison of $\phi-\tilde{\theta}$ distributions at polar radius $r=0.2$.}
    \label{Comp_u_02}
\end{figure}

\begin{figure}[h!]
    \centering
    \subfigure[Comparison of the real part for a polar radius of $r=0.4$.]{
    \centering
    \includegraphics[width=1\linewidth,height=0.22\textwidth]{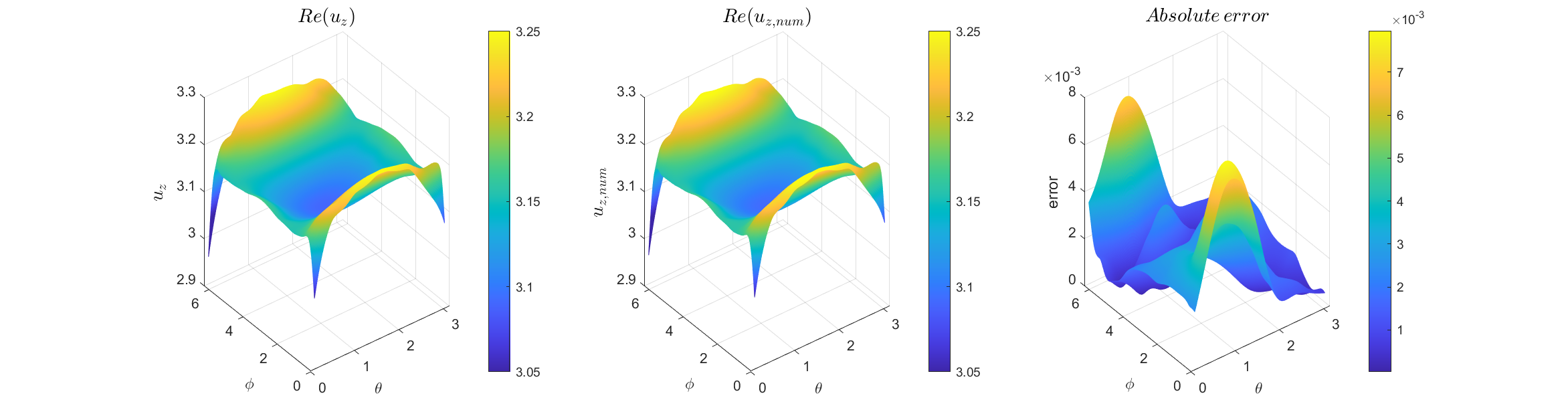}
    }
    \subfigure[Comparison of the imaginary part for a polar radius of $r=0.4$.]{
    \centering
    \includegraphics[width=1\linewidth,height=0.22\textwidth]{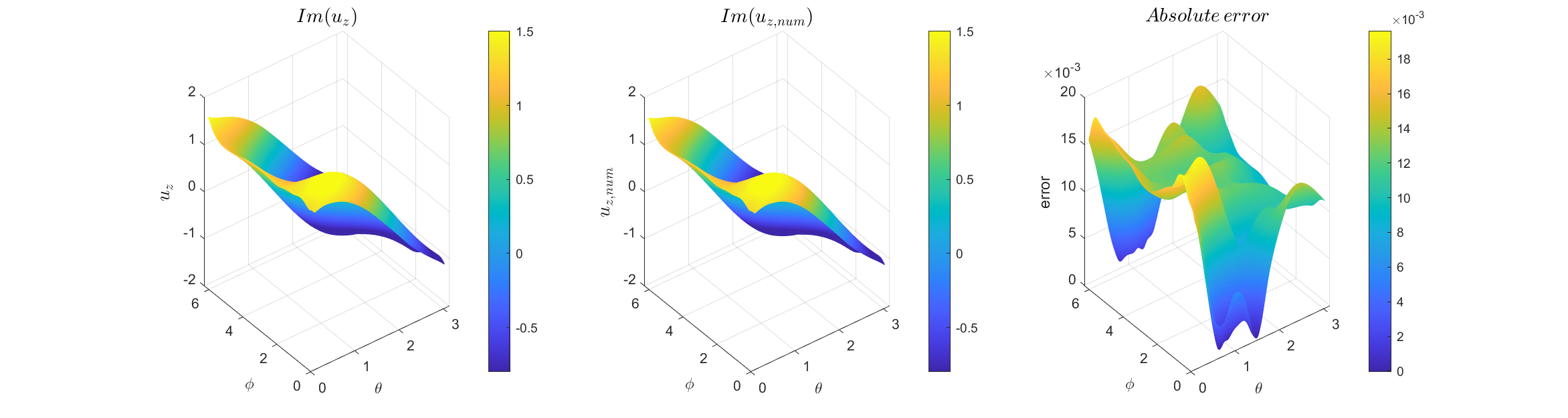}
    }
    \caption{Comparison of $\phi-\tilde{\theta}$ distributions at polar radius $r=0.4$.}
    \label{Comp_u_04}
\end{figure}

From Figure \ref{Comp_u_02} and Figure \ref{Comp_u_04}, it is clear to see that, for $n=200$, the real part errors $E_{\mathrm{pw}}^{\mathrm{re}}(x)$ remain around $10^{-3}$, while the imaginary part errors $E_{\mathrm{pw}}^{\mathrm{im}}(x)$ exhibit localized peaks reaching the order of $10^{-2}$.
Furthermore, we can also obtain the absolute error estimate in the full domain as
\begin{equation}\label{WZEq2}
    E_2=\left\Vert u_{z}(\cdot, \tilde{\theta}, \phi, \theta, \omega) - u_{z,num}( \cdot, \tilde{\theta}, \phi, \theta, \omega) \right\Vert_{\mathbb{L}^{2}(B({\bf 0},1))}  \simeq  0.051.
\end{equation}

\begin{figure}[h!]
    \centering
    \subfigure[Comparison of real part of the far-field pattern.]{
    \centering
    \includegraphics[width=1\linewidth,height=0.22\textwidth]{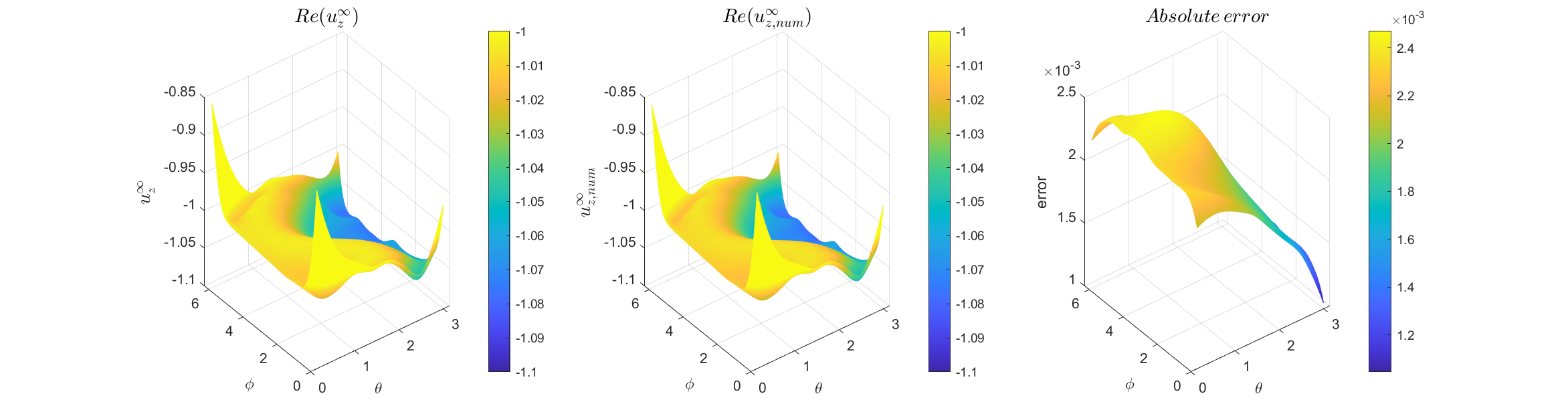}
    }
    \subfigure[Comparison of imaginary part of the far-field pattern.]{
    \centering
    \includegraphics[width=1\linewidth,height=0.22\textwidth]{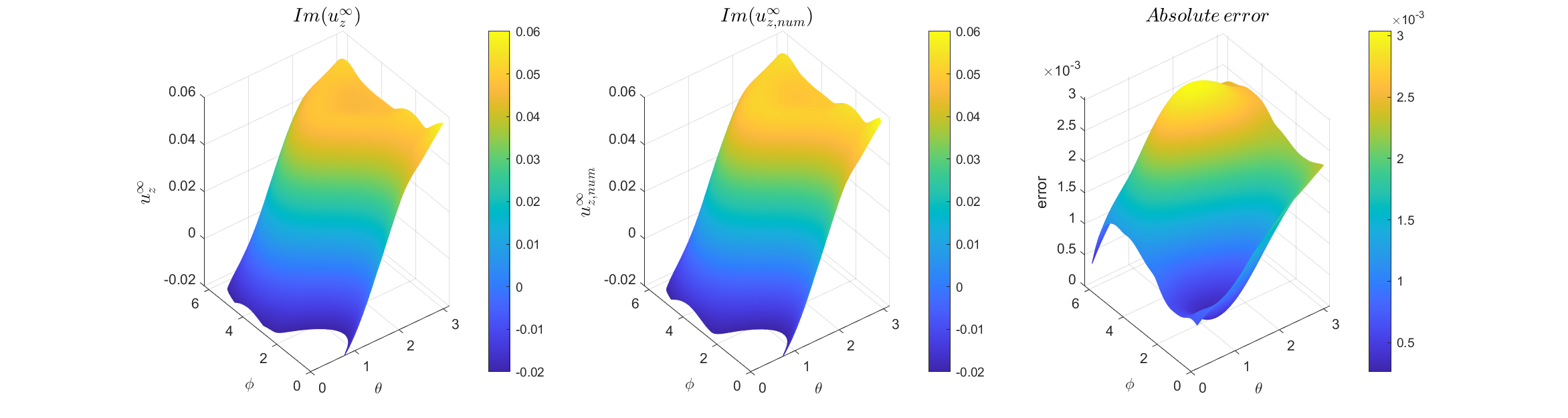}
    }
    \caption{Comparison of the far-field for $\hat x\in \mathbb{S}^2$.}
    \label{Compuinf}
\end{figure}

It is clear that our algorithm also has high accuracy in calculating the the perturbed medium case.
%by a droplet case.
Therefore we can use the coefficient $(\beta_{1}, \cdots, \beta_{n} )$ of \eqref{liu011} to compute $u^{\infty}_{z,num}(\hat{x}, \theta, \phi)$ to compare the far field $u^{\infty}_{z}(\hat{x}, \theta, \phi)$ generated by $ u_{z}(\cdot,\cdot,\cdot)$, which are shown in Figure \ref{Compuinf}.

\begin{remark}
    Although the technique used for both perturbed medium cases and unperturbed medium cases is identical, there is a slight discrepancy in the numerical solutions that are generated, see  \eqref{WZEq1} and  \eqref{WZEq2}. This is natural and due to the presence of the droplet $D_{z}$, with small size and high contrasting bulk coefficient, inside the domain $\Omega \equiv B({\bf 0},1)$.
\end{remark}

In this section, we validated the numerical solver by generating known solutions using analytic test functions with a fixed number of collocation points, $n = 200$. This test confirms the correctness of the solver itself but does not address the convergence order or robustness of the proposed scheme. To assess these, we must execute the solver for various values of $n$ and for the perturbation values of $\Theta(x,\theta,\omega)$ given by \eqref{FFJM}.

\section{Recovering the bulk modulus function $k_{0}(\cdot)$ inside $\Omega \, \equiv \, B({\bf 0},1)$}\label{SectionV}
The goal of this section is to numerically recover the bulk modulus function $k_{0}(\cdot)$ in local domain of $\Omega \, \equiv \,B({\bf 0},1)$, using the scheme described by \textbf{Algorithm \ref{AlgoScheme}} from the data
\begin{equation}\label{BF-Equa0354}
       \xi(z)  =  v^{\infty}(-  \theta, \theta, \omega_{1})  -  u^{\infty}_z(-  \theta, \theta, \omega_{1}),
    \end{equation}
see \eqref{KK}. It is important to note that in engineering situations, the above inversion input data $\xi(\cdot)$ are obtained through physical measurements for $\left( v^{\infty}(-  \theta, \theta, \omega_{1}),\, u^{\infty}_{z}(-  \theta, \theta, \omega_{1}) \right)$. To test our reconstruction algorithm, we use the far-field data $v^{\infty}(-\theta, \theta, \omega_{1})$ and $u^{\infty}_z(-  \theta, \theta, \omega_{1})$ generated from the simulation process. The effectiveness of the simulation process has been checked in \textbf{Section \ref{Sec-Num-Sim}}.  We assume the bulk modulus function $k_{0}(\cdot)$ to be identified is given by
\begin{equation}\label{k0GRE}
    k_{0,exact}(x) = \frac{2}{1+\left| x \right|^2} , \quad x \,\in B({\bf 0},1),
\end{equation}
for the simulation. Following the arguments in \textbf{Section \ref{SectionII}}, for an assumed known bulk modulus function, see for instance $(\ref{k0GRE})$, we can generate the far-field quantities $v^{\infty}(\cdot)$ and $u^{\infty}(\cdot)$, and consequently $\xi(z)$ as given in $(\ref{BF-Equa0354})$, for $z\in B({\bf 0},1)$. Specifically, as explained in \textbf{Section \ref{SectionII}}, the simulation of $u^\infty_z(-\theta,\theta,\omega_1)$ for $z\in \mathfrak{U}$ where $k_0(z)$ is recovered, requires to solve a 3-dimensional linear integral equation for every $z\in \mathfrak{U}$, which need large computational cost for large domain $\mathfrak{U}$. To check the numerical performances of  \textbf{Algorithm \ref{AlgoScheme}} with acceptable computational costs, we take  $\mathfrak{U}$ as a cube $\mathfrak{Q}:\equiv \mathfrak{Q}({\bf 0};0.5)\subset\subset B({\bf 0},1)$, which is centered at the origin with the side length  $0.5$. In this cube, we take $61$-points along each direction as observation points.
More precisely, the periodically distributed discrete points $x_{i,j,k} = \left(x^1_{i},x^2_{j},x^3_{k}\right) \in \mathfrak{Q}$ are specified by
\begin{equation*}
     x^1_{i}=-0.25+\frac{0.5(i-1)}{N-1},\quad x^2_{j}=-0.25+\frac{0.5(j-1)}{N-1}\quad \text{and} \quad x^3_{k}=-0.25+\frac{0.5(k-1)}{N-1}
\end{equation*}
for $ i, j, k=1,\cdots, N$ with $N=61$, and the corresponding noisy data $\xi^{\tau }(\cdot)$ are obtained by shifting  $\xi(\cdot)$ using a random variable $\mathcal{X}(\cdot)$ following a uniform distribution $\mathcal{U}\left( \cdot \right)$, i.e., $\mathcal{X}  \sim  \left( \dfrac{\tau }{0.25} \right)^{3}  \mathcal{U}\left( \mathfrak{Q} \right)$, with $\tau $ being a small positive number indicating the noisy level. For $\left\{ \tau _{i,j,k} := \mathcal{X}\left( x_{i,j,k} \right):\,  i, j, k=1,\cdots, N \right\}$, we synthesize the noisy data by
\begin{equation*}\label{noisydata}
    \xi^{\tau }\left(x_{i,j,k}\right)  := \xi\left(x_{i,j,k}\right) +  \tau _{i,j,k}  \xi\left(x_{i,j,k}\right) \quad \text{for} \quad  i, j, k=1,\cdots, N \quad \text{and} \quad 0 < \tau _{i,j,k} \ll 1,
\end{equation*}
which fulfills the inequality
\begin{equation*}
    \left\Vert \xi^{\tau }(\cdot)-\xi(\cdot) \right\Vert_{\mathbb{L}^{\infty}\left(\mathfrak{Q} \right)}  \le  \tau \left\Vert \xi(\cdot)\right\Vert_{\mathbb{L}^{\infty}\left(\mathfrak{Q} \right)}.
\end{equation*}

Although we only check our reconstruction algorithm in an interior local domain $\mathfrak{Q}=\mathfrak{Q}({\bf 0}, 0.5)\subset B({\bf 0},1)$ due to the restrictions on the computational cost for yielding the inversion input by simulation process, our algorithm can be applied to any interior domain of $\Omega$, if we do not consider the simulation cost for yielding $u^{\infty}_{z}(\cdot,\cdot,\cdot)$, or we assume that $u^{\infty}_{z}(\cdot,\cdot,\cdot)$ are given directly in $\Omega \equiv B({\bf 0},1)$.
\medskip
\newline
The crucial step for numerically computing $k_0(z)$ by $(\ref{RF})$ from noisy inversion input is the computational scheme for $\Delta \xi^\tau (z), \nabla \xi^\tau (z)$ with noisy level $\tau >0$. By our mollification
scheme \eqref{grexi} and \eqref{lapxi}, we are essentially to compute  $\Delta \xi ^{\delta,\tau }(z)$  and $\nabla\xi^{\delta,\tau }(z)$ with specified $\delta>0$ and noisy level $\tau $.
\medskip
\newline
However, we only have the  noisy values $\xi^\tau $ at $61$ discrete points, due to the restriction on the computational cost in simulating these inversion input by solving 3-dimensional integral equations for each point $z$. These values are not enough for computing the derivatives accurately from \eqref{grexi} and \eqref{lapxi}, since the integrations in \eqref{grexi} and \eqref{lapxi} need very fine grids of $\xi^\tau $ for yielding good approximation. To overcome this limitation, we refine the computational grids by interpolation process using $\{\xi^\tau (x_{i,j,k}): i,j,k=1,\cdots, 61\}$ to yield a smooth function $\tilde\xi^\tau $, which is taken as the approximation to $\xi^\tau (x)$. This process is finished by the standard interpolation function \textbf{griddedInterpolant} in Matlab with cubic spline function as base functions. Then, we have the values of
$\tilde \xi^\tau $  at fine grids $\{\tilde x_{i,j,k}:\;i,j,k=1,\cdots,201\}\subset \mathfrak{Q}$. Thus, similarly to $(\ref{RF})$, we have the following reconstruction formula
\begin{equation}\label{XC-Equa0850dis}
\frac{1}{k_{0,num}(\tilde x_{i,j,k})} \approx - \frac{1}{\omega_{1}^2} \left( \frac{\Delta\tilde \xi^{\delta,\tau }(\tilde x_{i,j,k})}{2\tilde \xi^\tau (\tilde x_{i,j,k})}-\frac{|\nabla\tilde \xi^{\delta,\tau }(\tilde x_{i,j,k})|^2}{4\tilde \xi^\tau (\tilde x_{i,j,k})^2}\right) , \quad \tilde x_{i,j,k}\in \mathfrak{Q}.
\end{equation}

In the numerical experiments, we choose the mollification radius proportional to the grid size by $\delta=c\frac{0.5}{N}$, to ensure the stable numerical quadrature, where $c=1,2,\cdots$. A large value of $\tau$ necessitates a large $c$ to suppress noise amplification in $\Delta\tilde \xi^{\delta,\tau }$ and $\nabla\tilde \xi^{\delta,\tau }$. Besides, since we have the values of $\tilde\xi^\tau (\cdot)$ at very fine grids $\{\tilde x_{i,j,k}:\;i,j,k=1,\cdots,201\}$, the right-hand side of $(\ref{XC-Equa0850dis})$ can be computed accurately, as mentioned in $(\ref{grexi})$ and $(\ref{lapxi})$ .
\medskip
\newline
To indicate the numerical reconstruction performance of the \textbf{Algorithm \ref{AlgoScheme}} quantitatively,
we introduce the global relative error (GRE) and the point-wise relative error (PRE) in $\mathfrak{Q}$ related to the  {\it a-prior}  known bulk modulus function $k_{0,\text{exact}}(\cdot)$, given by $(\ref{k0GRE})$,
\begin{equation*}\label{Def-GRE}
    GRE :=  \frac{\left\Vert k_{0,exact}  -  k_{0,num} \right\Vert_{\ell^2\left(\mathfrak{Q}\right)}}{\left\Vert k_{0,exact}\right\Vert_{\ell^2\left(\mathfrak{Q}\right)}}  :=  \frac{\left( \displaystyle\sum_{i=1}^{201} \sum_{j=1}^{201} \sum_{k=1}^{201} \left\vert k_{0,exact}(\tilde x_{i,j,k})  -  k_{0,num}(\tilde x_{i,j,k}) \right\vert^{2} \right)^{\frac{1}{2}}}{\left( \displaystyle\sum_{i=1}^{201} \sum_{j=1}^{201} \sum_{k=1}^{201} \left\vert k_{0,exact}(\tilde x_{i,j,k}) \right\vert^{2} \right)^{\frac{1}{2}}},
\end{equation*}
and
\begin{equation*}\label{Def-PRE}
    PRE(\tilde x_{i,j,k}) :=  \frac{\left| k_{0,exact}(\tilde x_{i,j,k}) - k_{0,num}(\tilde x_{i,j,k})\right|}{\left| k_{0,exact}(\tilde x_{i,j,k})\right|}, \quad
    \tilde x_{i,j,k}\in \mathfrak{Q},
\end{equation*}
respectively,
where $ k_{0,num}(\cdot)$ is the numerical solution obtained from inverting the data \eqref{BF-Equa0354}. Since the mollification scheme treating the numerical differentiations for noisy data is of only average convergence by \eqref{VUPR}, we can also expect the good approximation under the convergence norm, instead of the point-wise approximation.

\begin{figure}[H]
    \centering
    \includegraphics[width=1\linewidth,height=0.22\textwidth]{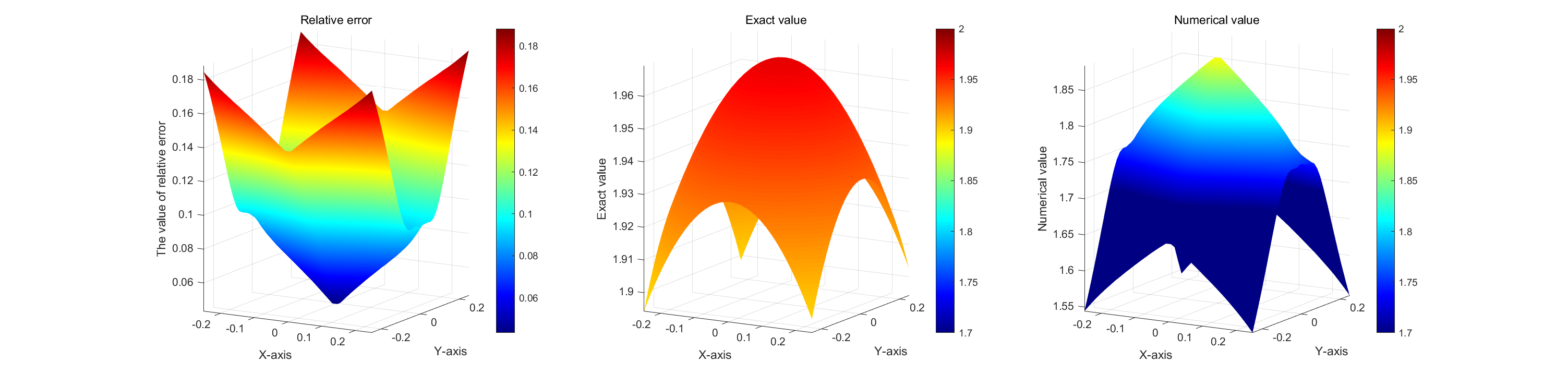}
    \caption{\scriptsize  The recovery of $k_{0,exact}(\cdot)$ using exact inversion input with $\tau  = 0$ on $x_3=0.125$.}
    \label{k0n0}
\end{figure}

\begin{figure}[H]
    \centering
    \includegraphics[width=1\linewidth, height=0.22\textwidth]{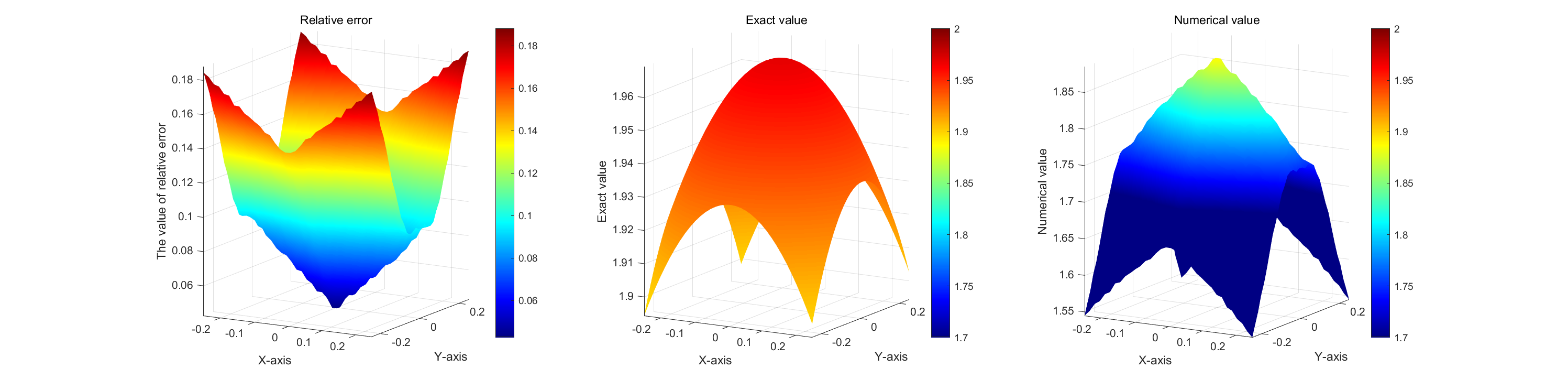}
    \caption{\scriptsize  The recovery of $k_{0,exact}(\cdot)$ using noisy inversion input with $\tau  = 0.01$ on $x_3=0.125$.}
    \label{k0n01}
\end{figure}

\begin{figure}[H]
    \centering
    \includegraphics[width=1\linewidth, height=0.22\textwidth]{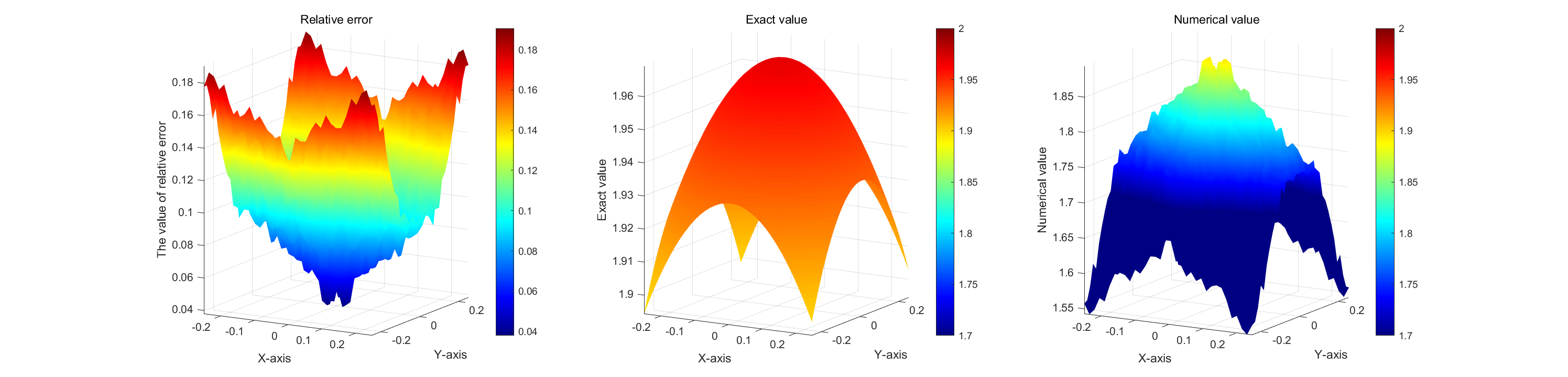}
    \caption{\scriptsize The recovery of $k_{0,exact}(\cdot)$ using noisy inversion input with $\tau  = 0.05$ on $x_3=0.125$.}
    \label{k0n05}
\end{figure}

\begin{figure}[H]
    \centering
    \includegraphics[width=1\linewidth, height=0.22\textwidth]{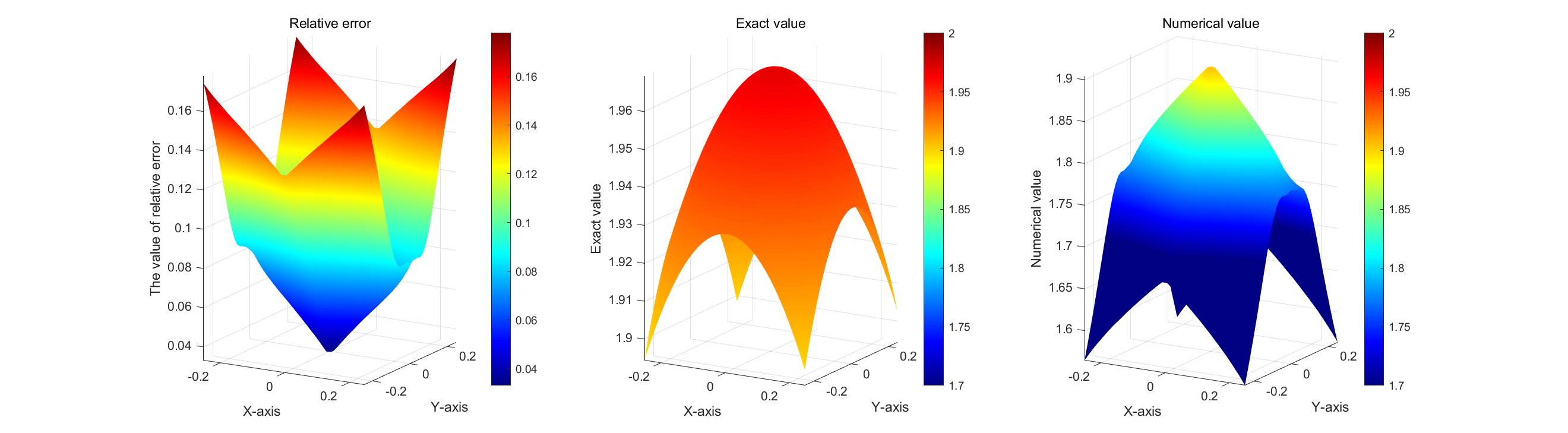}
    \caption{\scriptsize  The recovery of $k_{0,exact}(\cdot)$ using exact inversion input with $\tau  = 0$ on $x_2 =-0.125$.}
    \label{k0x2n0}
\end{figure}

\begin{figure}[H]
    \centering
    \includegraphics[width=1\linewidth, height=0.22\textwidth]{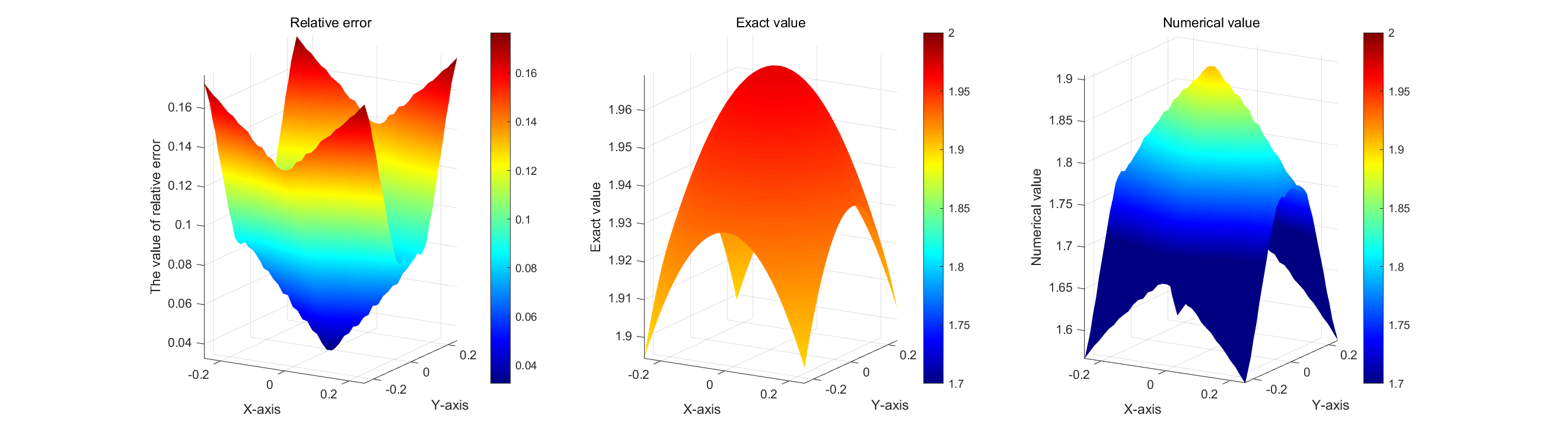}
    \caption{\scriptsize  The recovery of $k_{0,exact}(\cdot)$ using noisy inversion input with $\tau  = 0.01$ on $x_2 =-0.125$.}
    \label{k0x2n01}
\end{figure}

\begin{figure}[H]
    \centering
    \includegraphics[width=1\linewidth, height=0.22\textwidth]{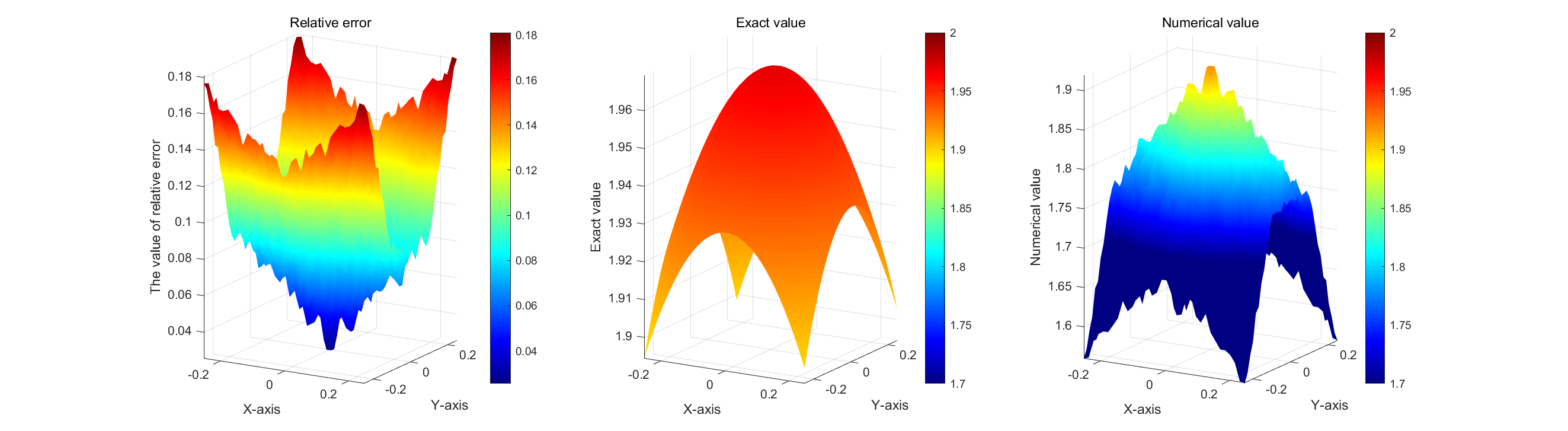}
    \caption{\scriptsize  The recovery of $k_{0,exact}(\cdot)$ using noisy inversion input with $\tau  = 0.05$ on $x_2 =-0.125$.}
    \label{k0x2n05}
\end{figure}
Noticing that $k_0(\cdot)$ is a functional with 3-dimensional argument, for the ease of presentation, we show the reconstructions together with exact  $k_{0,exact}(\cdot)$ only on the cross planes $x_3 = 0.125$ and $x_2 = -0.125$ for noise levels $\tau =0,0.01,0.05$.
Figure \ref{k0n0} to Figure \ref{k0x2n05}, the leftmost subplot of each figure displays the point-wise relative error, the central subplot shows the distribution of the exact bulk modulus function $k_{0,exact}(\cdot)$, and the rightmost subplot presents the numerical reconstruction $k_{0,num}(\cdot)$ obtained using noisy data.

We observe that in the absence of artificially added noise, the exact solution $k_{0,exact}(\cdot)$ and the numerical solution $k_{0,num}(\cdot)$ exhibit identical distributions, albeit with slight variations in magnitude.
On the cross-plane $x_3 = 0.125$, the maximum point-wise relative error remains below $0.18$. The GRE in the absence of noise is $0.1081$, these results demonstrate the high accuracy and feasibility of our method. When noise at a level of $0.01$ is introduced, the numerical solution $k_{0,num}(\cdot)$ still maintains a similar distribution, and the point-wise relative error stays below $0.19$, at the same time the global relative error is $0.1126$, albeit with increased oscillations. At a noise level of $0.05$, the reconstruction exhibits  oscillation. Nevertheless, in this case the global relative error is $0.1183$.
Additionally, on the cross-plane $x_2 = -0.125$, the maximum point-wise relative error remains below $0.17$, and the GRE is $0.1068$. This case shows good accuracy when no noise is added. For $\tau=0.01$, the maximum point-wise relative error is below $0.19$, and the GRE is $0.1083$. At a noise level of $0.05$, the reconstruction has more visible oscillations. However, the GRE is still only $0.1105$.

\begin{table}[h]
\centering
\caption{GRE for different noise levels on planes $x_3 = 0.125$ and $x_2=-0.125$. }
\begin{tabular}{c|c|c|c|c|c}
\hline
Noise level $\tau$ & $0$ & $0.01$ & $0.05$ & $0.1$ & $0.15$ \\
\hline
GRE on cross-plane $x_3 = 0.125$ & $0.1081$ & $0.1126$ & $0.1183$ & $0.1278$ & $0.1465$\\
\hline
GRE on cross-plane $x_2 = -0.125$ & $0.1068$ & $0.1083$ & $0.1105$ & $0.1189$ & $0.1401$\\
\hline
\end{tabular}
\label{grex3x2}
\end{table}

More information about the global relative errors for different noise levels on the cross-plane $x_3 = 0.125$, and the cross-plane $x_2 = -0.125$ is shown in Table \ref{grex3x2} quantitatively. These results show that the reconstruction is not very sensitive to the noise level. Therefore, the method is stable in our numerical experiments. Finally, our numerical experiments show that the reconstruction is stable for $\delta$ within a practical interval, which means the robustness of the mollification regularization in our setting.

\vskip 0.5cm

{\bf Acknowledgment.} The authors thank the referees very much for their valuable comments and suggestions, which make this paper much improved.

\vskip 0.5cm

\end{document}